\documentclass[12pt,reqno]{amsart}
\usepackage[utf8]{inputenc} 
\usepackage[T1]{fontenc}
\usepackage{xifthen,csquotes,amsmath}
\usepackage{extarrows}
\usepackage{multirow,bigints,amsmath}
\usepackage{subfig}
\let\Horig\H

\usepackage{pkgfile}

\DeclareMathOperator\arctanh{arctanh}

\begin{document}

\title[Statistics of the two-star ERGM]{Statistics of the two-star ERGM}

\author[Sumit Mukherjee and Yuanzhe Xu]{Sumit Mukherjee \& Yuanzhe Xu}
\thanks{Research partially supported by NSF Grant DMS-1712037}

\address{Department of Statistics, Columbia University
\newline\indent 1255 Amsterdam Avenue, New York, NY 10027}

\date{\today}



\maketitle

\begin{abstract}
%

In this paper, we explore the two-star Exponential Random Graph Model, which is a two parameter exponential family on the space of simple labeled graphs. We introduce auxiliary variables to express the two-star model as a mixture of the $\beta$ model on networks. Using this representation, we study asymptotic distribution of the number of edges, and the sampling variance of the degrees. In particular, the limiting distribution for the number of edges has similar phase transition behavior to that of the magnetization in the Curie-Weiss Ising model of Statistical Physics.  Using this, we show existence of consistent estimates for both parameters in all parameter domains. Finally, we prove that the centered partial sum of degrees converges as a process to a Brownian bridge in all parameter domains, irrespective of the phase transition. 
\end{abstract}
\section{Introduction}

Inference on graphs/networks is a topic of considerable recent interest in Statistics and Machine Learning. Both parametric and non-parametric models have been introduced to model graphs. In the parametric setting, perhaps the simplest model is the celebrated Erd\Horig{o}s-R\'enyi model, where all the edges are independent, and there is only one parameter in the model. However, this model is too simplistic to be able to capture real life networks. Note that an Erd\Horig{o}s-R\'enyi model can be expressed as an exponential family, with the number of edges  as a sufficient statistic. As a first step towards modeling dependence between edges, it is natural to consider parametric models where there are more than one sufficient statistic. Typical sufficient statistics on graphs of interest include higher order subgraph counts, such as number of stars, number of cycles, number of cliques, and so on. 
 This motivation led to the introduction and study of exponential families on the space of graphs with finitely many sufficient statistics. We will refer to these class of models as 
 Exponential Random Graph Models. For the sake of convenience the abbreviation ERGM will henceforth be used to refer to Exponential Random Graph Models. 
 ERGMS first appeared in Social Sciences (c.f.~\cite{anderson1999p,frank1986markov, holland1981exponential,robins2007introduction,wasserman1994social,wasserman1996logit} and references there-in), and since then have received a lot of attention in Probability (c.f.~\cite{ganguly2019sub,gotze2019concentration,radin2013phase} and references there-in), Statistics (c.f.~\cite{chatterjee2013estimating,schweinberger2020concentration,shalizi2013consistency} and references there-in) and Statistical Physics (\cite{demuse2018phase,park2004solution,park2004statistical} and references there-in). 
 \\
 
 Among the class of ERGMs, perhaps the simplest is the two-star model, first studied in \cite{park2004solution}. This ERGM has exactly two sufficient statistics, the number of edges, and the number of two stars.   
 By a two-star, we mean a path of length $2$, which has $3$ vertices and $2$ edges. The advantage of working with two stars is that the number of two stars can be expressed very conveniently as a function of the degrees of the graph, as we will see below.
 \\
 
 We begin by introducing the two-star ERGM.
\begin{defn}
For a positive integer $n$, let $\mathcal{G}_n$ denote the space of all simple graphs with vertices labeled $[n]:=\{1,2,\ldots,n\}$. Since a simple graph is uniquely identified by its adjacency matrix, without loss of generality we can take $\mathcal {G}_n$ to also denote the set of all symmetric $n\times n$ matrices, with $0$ on the diagonal elements and $\{0,1\}$ on the off-diagonal elements. By slightly abusing the notation, we use $G$ to denote both a graph and its $n\times n$ adjacency matrix $(G_{ij})_{1\le i,j\le n}$, defined by 
\begin{align*}
G_{ij}=\left\{
\begin{array}{rcl}
1&     &\text{If an edge is present  between vertices $i$ and $j$ in $G$}\\
0&     &  \text{Otherwise}
\end{array} \right.
\end{align*} 
Set $G_{ii}=0$ by convention. Let $E(G):=\sum_{i<j}G_{ij}$ denote the number of edges in $G$, and let 
\begin{align*}
T(G):=\sum_{i=1}^n\sum_{j<k}G_{ij}G_{ik}
\end{align*}
denote the number of two-stars in $G$. A simple calculation shows that the number of two stars can be written as
\begin{align*}
T(G)=\sum_{i=1}^n{d_i(G)\choose 2},
\end{align*}
 where $(d_1(G),\ldots,d_n(G))$ is the labeled degree sequence of the graph $G$, defined by $d_i(G):=\sum_{j=1}^nG_{ij}$. Indeed, this is because given any vertex $i$ of degree $d_i(G)$, there are ${d_i(G)\choose 2}$ two stars with $i$ as their central vertex.
 \\
 
 Given parameters $\omega_1>0$ and $\omega_2\in \R$, the two-star ERGM is defined by the following probability mass function on $\mathcal{G}_n$:
\begin{align}\label{eq:edge_star}
\P_{n}(G=g):=\frac{1}{{Z}_n(\omega_1,\omega_2)} \text{exp}\Big\{\Big(\omega_2+\frac{\omega_1}{n-1}\Big)E(g)+\frac{\omega_1}{n-1}T(g)\Big\},
\end{align} where ${Z}_n(\omega_1,\omega_2)$ is the normalizing constant.
\end{defn}
Note that if $\omega_1=0$, the model reduces to an Erd\Horig{o}s-R\'enyi model with parameter $p=\frac{e^{\omega_2}}{1+e^{\omega_2}}.$ The regime $\omega_1>0$ corresponds to the so called \enquote{Ferromagnetic regime} of Statistical Physics, which encourages more two stars in the graph than an Erd\Horig{o}s-R\'enyi graph with parameter $p$ as above. 
One of the main difficulties of analyzing this model (and ERGMs in general) is that the normalizing constant $Z_n(\omega_1,\omega_2)$ is not available in closed form. Explicit computation of the normalizing constant is computationally prohibitive. One way out is to resort to MCMC, but mixing rates for ERGMs depend crucially on the parameter values as shown in \cite[Theorem 5,6]{bhamidi2008mixing}, and can take time which is exponential in $n$ to mix. In \cite{chatterjee2013estimating} the authors study ERGMS using a large deviation approach. In particular they show that the two star ERGM is \enquote{close} to a mixture of Erd\Horig{o}s-R\'enyi random graphs, for all $(\omega_1,\omega_2)\in \R^2$ (c.f.~\cite[Theorem 6.4]{chatterjee2013estimating}). Since an Erd\Horig{o}s-R\'enyi random graph has a single parameter, this suggests that consistent estimation of both the parameters (as the size of the graph grows) is not possible in the two star ERGM.
Nevertheless, as we show below in Corollary \ref{cor:est}, it is possible to estimate both the parameters in a consistent manner, though there is a loss of efficiency when trying to estimate two parameters instead of one. The next section summarizes our main results.

\subsection{Main results}\label{Results}
Throughout the paper, we work in a slightly different parametrization, given by
\begin{align}\label{eq:reparameter}
\theta:=\frac{\omega_1}{4}>0, \beta:=\frac{\omega_1+\omega_2}{2}\in \R.
\end{align}
As frequently happens for such models, the two star model undergoes a phase transition, and its behavior is qualitatively different in different parts of the parameter regime. The following lemma introduces these different parameter domains. The proof of this lemma follows from straightforward calculus, and is deferred to the appendix (section \ref{sec:four1}). 
\begin{lem}\label{lem:min}
Setting
\begin{align}\label{eq:central}
q(x)=\theta x^2-\log\cosh(2\theta x+\beta),
\end{align}
the following hold:
\begin{enumerate}
\item[(a)]
If either $\theta>0,\beta\ne 0$ or $\theta\in (0,1/2), \beta=0$, the function $q(.)$ has a unique global minimizer at $t$, where $t$ is the unique root of the equation $x=\tanh(2\theta x+\beta)$ which has the same sign as that of $\beta$. Further we have $q''(t)=2\theta[1-2\theta(1-t^2)]>0$. 

\item[(b)]
If $\theta>1/2,\beta=0$, the function $q(.)$ has two global minimizers at $\pm t$, where $t$ is the unique positive root of the equation $x=\tanh(2\theta x)$. Further, we have $q''(\pm t)=2\theta[1-2\theta(1-t^2)]>0$.
\item[(c)]
If $\theta=1/2, \beta=0$, the function $q(.)$ has a unique global minimizer at $t=0$, and $q''(0)=0$. 

\end{enumerate}
\end{lem}
\begin{defn}
Let $t$ be as defined in Lemma \ref{lem:min}, and note that $t$ depends on $(\theta,\beta)$, which we suppress for ease of notation. Also let 
\begin{eqnarray*}
&\Theta_{11}:=\{\theta\in (0,1/2), \beta=0\},\quad  &\Theta_{12}:= \{\theta>0, \beta\ne 0\}\\
& \Theta_2:=\{\theta>1/2, \beta=0\},\quad &\Theta_3:=\{\theta=1/2,\beta=0\},
 \end{eqnarray*}
 and set $\Theta_1:=\Theta_{11}\cup \Theta_{12}$. We will refer to the three regimes $\Theta_1, \Theta_2, \Theta_3$ as uniqueness regime, non uniqueness regime, and critical regime respectively, the reason for this nomenclature follows from Lemma \ref{lem:min}.
\end{defn}

Our first main result now gives the asymptotic distribution for the number of edges in all the three domains $\{\Theta_1, \Theta_2, \Theta_3\}$.
\begin{thm}\label{mean1}
Suppose $G$ is a random graph from the two star model in \eqref{eq:edge_star}. 
\begin{enumerate}
\item[(a)]
If $(\theta,\beta)\in \Theta_1$, we have
\begin{align}\label{Asymptotic1}
n\Big(\frac{2E(G)}{n^2}-p\Big)\stackrel{d}{\longrightarrow} N\Big(-\mu,\sigma^2\Big),
\end{align}
where $p=\frac{1+t}{2}$, $\mu:=\frac{\theta t(1-t^2)}{[1-\theta(1-t^2)][1-2\theta(1-t^2)]}$, and $\sigma^2:=\frac{1-t^2}{2-4\theta(1-t^2)}$.
%
%
%
%
%
\item[(b)]
If $(\theta,\beta)\in \Theta_2$, then we have
\begin{align}\label{Asymptotic2}
&n\Big(\frac{2E(G)}{n^2}-p\Big)\Big|E(G)> \frac{n^2}{4}\Big) \stackrel{d}{\longrightarrow}N(\mu,\sigma^2),
\end{align}
in which $p, \mu, \sigma^2$ have the same formulas as above. 


\item[(c)]
If $(\theta,\beta)\in \Theta_3$, then we have
\begin{align}\label{Asymptotic3}
2\sqrt{n}\Big(\frac{2E(G)}{n^2}-\frac{1}{2}\Big)\stackrel{d}{\longrightarrow}\zeta,
\end{align}
where $\zeta$ is a random variable on $\R$ with density proportional to $e^{-\zeta^2/2-\zeta^4/24}$.
\end{enumerate}
\end{thm}

\begin{remark}
Prior to our work, limiting distribution for $E(G)$ was not understood for any ERGM. See however the recent work of \cite{ganguly2019sub}, which studies asymptotic distribution of the sum of a small number of disjoint edges, in the high temperature regime (i.e. $\theta$ small). Also related to our work is the asymptotics of the magnetization/sum of spins in Ising models on dense regular graphs. As explained below in section \ref{sec:ising}, the two star ERGM can be thought of as an Ising model on a $d_N$ regular graph with $N:={n\choose 2}$ vertices and degree $d_N=2(n-2)$ (so that $d_N\propto \sqrt{N}$). Further, the number of edges is a linear function of the magnetization. Very recently in \cite{deb2020fluctuations} it was shown that the magnetization/sum of spins in a Ising model on a regular graphs with degree $d_N\gg\sqrt{N}$ is universal, and is the same as obtained for the Curie-Weiss model ($d_N=N-1$) in \cite{ellis1978statistics}. As our results demonstrate, universality breaks at the threshold $d_N\propto \sqrt{N}$, as the distribution above does not match that of the Curie-Weiss model in the domains $\Theta_{12}\cup \Theta_2\cup \Theta_3$. Only in the domain $\Theta_{11}$ the limiting distribution of the magnetization in the two star model matches that of the Curie-Weiss model. The techniques employed in this draft are very different from the techniques of both \cite{ganguly2019sub} and \cite{deb2020fluctuations}.
\end{remark}
Our second result studies the fluctuations of the empirical variance of the degrees. 
\begin{thm}\label{sd}
Suppose $G$ is a random graph from the two star model in \eqref{eq:edge_star}. For all $\theta>0, \beta\in \R$ we have
\begin{align}\label{sd1}
\sqrt{n}\left[\frac{4}{n^2}\sum\limits_{i=1}^n(d_i(G)-\bar{d}(G))^2-\tau\right]\stackrel{d}{\longrightarrow}N(0,2\tau^2)
\end{align}
where $\bar{d}(G):=\frac{\sum_{i=1}^{n}d_{i}(G)}{n}$ and $\tau:=\frac{1-t^2}{1-\theta(1-t^2)}$.

\end{thm}
\begin{remark}
In particular the empirical variance of the degrees converge in probability to $\tau$, which is continuous but not differentiable at $\theta=1/2$, when $\beta=0$ is kept fixed. This phenomenon was also observed in \cite[Fig 2]{park2004solution}. 
\end{remark}
As an application of the two theorems above, we provide consistent estimators of the parameters $(\theta,\beta)$.
\begin{cor}\label{cor:est}

Suppose $G$ is a random graph from the two star model in \eqref{eq:edge_star}, with $(\theta,\beta)\in \Theta$.

\begin{enumerate}
\item[(a)]
If $(\theta,\beta)$ are both unknown, then there exists a $\sqrt{n}$ consistent estimator for $(\theta,\beta)$, i.e. $\sqrt{n}(\hat{\theta}-\theta,\hat{\beta}-\beta)=O_P(1)$.

\item[(b)]
If $\theta$ is known, then there is a $n$ consistent estimator for $\beta$, i.e. $n(\hat{\beta}-\beta)=O_P(1)$.
\end{enumerate}

\end{cor}

%
\begin{remark}
The above corollary shows that there is a loss of efficiency when we are trying to estimate both parameters, as opposed to estimating just one parameter. Joint estimation of parameters in general Ising models has been studied in \cite{ghosal2020joint}, where the authors give a general upper bound on the rate of consistency of pseudo-likelihood (see \cite[Theorem 1.2]{ghosal2020joint}). Using their result for a $d_N$ regular graph on $N$ vertices, one concludes that (an upper bound to) the rate of estimation error the pseudo-likelihood estimator is $\frac{d_N}{\sqrt{N}}$. Thus one can consistently estimate $(\theta,\beta)$ on an Ising model on a sequence of $d_N$ regular graph, if $d_N\ll \sqrt{N}$. However, in this case we have $d_N\propto \sqrt{N}$, and so consistency of the bivariate pseudo-likelihood estimator does not follow from \cite{ghosal2020joint}. It is unclear whether the pseudo-likelihood estimator is consistent in this case. On the other hand, the above corollary gives explicit consistent estimator for both parameters.
\end{remark}
%
Our final result shows that the partial sums of the (centered) degree distribution converges as a process in $\mathcal{C}[0,1]$ to a Brownian bridge under proper scaling, in all the three parameter domains. 
\begin{thm}\label{bm}
Suppose $G$ is a random graph from the two star model in \eqref{eq:edge_star}. Let $W_n(.)\in\mathcal{C}[0,1]$ be the linear interpolation of the points $\{(\frac{i}{n},\frac{S_i({\bf d})}{n-1}),i\in [n]\}$, where $S_i({\bf d}):=\sum\limits_{j=1}^i(d_j(G)-\bar{d}(G))$. Then \[W_n(.)\stackrel{d}{\longrightarrow}\sqrt{\tau}\{W(.)\},\] where $W(.)\in \mathcal{C}[0,1]$ is a Brownian bridge. 
\end{thm}
This demonstrates that irrespective of the phase transitions, there is significant Gaussian behavior in the model, which is captured in terms of contrasts. Similar Gaussian fluctuations were obtained in \cite{papangelou1989gaussian} for the Curie-Weiss model at criticality.

\subsection{ {Auxiliary Variables}}\label{sec:ising}

The main technique for proving the results of this paper is a representation of the two star model as a mixture of $\beta$ models by introducing auxiliary variables, introduced below.  We note that introducing auxiliary variables have been proved to be successful in rigorously analyzing the Curie-Weiss model (\cite{comets1991asymptotics,mukherjee2018global}), and have also been used in \cite{park2004solution} to study (non-rigorously) the two star model. Before introducing the auxiliary variable, we first transform the edge variables to $\{-1,1\}$ instead of $\{0,1\}$, and show that the transformed variables is a sample from an Ising model on an appropriate graph.
\\


Transform the edge variables from $\{0,1\}$ to $\{-1,1\}$ by setting $Y_{ij}:=2G_{ij}-1$ for $i\ne j$, and set $Y_{ii}:=0$ as convention. Via this transformation, the Hamiltonian for the matrix $Y:=(Y_{ij})_{1\le i,j\le n}$ (up to additive constants) is given by
\begin{align*}
&\frac{\omega_1}{4(n-1)}T(Y)+\frac{\omega_1+\omega_2}{2}E(Y)=\frac{\theta}{n-1}T(Y)+\beta E(Y),
\end{align*}
in which $\theta=\frac{1}{4}\omega_1$, and $\beta=\frac{1}{2}(\omega_1+\omega_2)\in \R$ as in Lemma \ref{lem:min}, and 
\begin{align*}
T(Y):=\sum\limits_{i=1}^n\sum\limits_{j<k}Y_{ij}Y_{ik},\quad E(Y)=\sum\limits_{i<j}Y_{ij}.
\end{align*} 
Thus the model $\P_n$ defined in \eqref{eq:edge_star} is an Ising model in the transformed variable $Y$ on the graph $\widetilde{G}_n$ which is the line graph of the complete graph $K_n$. More precisely, $\widetilde{G}_n$ has $\mathcal{E}:=\{(i,j)|1\leq i<j\leq n\}$ as its vertex set, and two distinct vertices $e=(i,j)$ and $f=(k,l)$ are connected iff $\{i,j\}\bigcap\{k,l\}\neq\emptyset$, i.e. $i=k$ or $i=l$ or $j=k$ or $j=l$. Thus $\widetilde{G}_n$ is a regular graph on ${n\choose 2}$ vertices, with degree $2(n-2)$.  Setting 
\begin{align}\label{eq:K}
k_i(Y):=\sum\limits_{j=1}^nY_{ij}=2d_i(G)-(n-1),
\end{align} the p.m.f. of $Y$ can be written as 
\begin{align}\label{Y:edge_star}
\P_{n}(Y=y)=\frac{1}{\widetilde{Z}_n(\theta,\beta)}\exp\left\{\frac{\theta}{2(n-1)}\sum\limits_{i=1}^{n}k_{i}(y)^{2}+\frac{\beta}{2}\sum\limits_{i=1}^{n}k_{i}(y)\right\}
\end{align}
Let $\phi=(\phi_1,\ldots,\phi_n)$ be a random vector in $\R^n$ defined by 
\begin{align}\label{eq:phi_z}
\phi_i=\frac{k_{i}(Y)}{n-1}+\frac{W_i}{\sqrt{(n-1)\theta}},  
\end{align}
where $(W_1,\ldots,W_n)\stackrel{i.i.d.}{\sim}N(0,1)$ are independent of the $Y$. 
The following proposition computes the distribution of $(Y|\phi)$, and the marginal density of $\phi$. The proof of this Proposition is deferred to the appendix (section \ref{sec:four1}).

\begin{ppn}\label{thm:bayes}
Suppose $Y$ is an observation from the p.m.f.~in \eqref{eq:edge_star}.
\begin{enumerate}
\item[(a)]
Given $\phi$, the random variables $\{Y_{ij}\}_{1\le i<j\le n}$ are mutually independent, with 
\begin{align}\label{y|phi}
\P_n(Y_{ij}=1|\phi)=\frac{e^{\theta(\phi_{i}+\phi_{j})+\beta}}{e^{\theta(\phi_{i}+\phi_{j})+\beta}+e^{-\theta(\phi_{i}+\phi_{j})-\beta}}    
\end{align}
\item[(b)]
The marginal density of $\phi$  has a density on $\R^n$ which is proportional to $e^{-f_n(\phi)}$, where 
$f_n(\phi):=\sum_{i<j} p(\phi_i,\phi_j)$ with
\begin{align}\label{eq:domm}
p(x,y)=\frac{\theta}{2}\Big(x^2+y^2\Big)-\log\cosh\Big[\theta(x+y)+\beta\Big]=\frac{\theta}{4}(x-y)^2+q\Big(\frac{x+y}{2}\Big),
\end{align}
where $q(x)=\theta x^2-\log\cosh(2\theta x+\beta)$ as in Lemma \ref{lem:min}.
\end{enumerate}

\end{ppn}


\begin{remark}
MCMC using auxiliary random variables is a common technique in simulations (\cite{andersen2007hit, edwards1988generalization, swendsen1987nonuniversal}). Using Proposition \ref{thm:bayes}, it follows that the conditional distribution of the graph $G$ given the vector $\phi$ is the $\beta-$model, which has received considerable attention in Statistics \cite{blitzstein2011sequential,chatterjee2019estimation,chatterjee2011random,rinaldo2013maximum} and references there-in). Thus the two-star model \eqref{eq:edge_star} can be expressed as a mixture of $\beta-$models with random weights.  Since both the conditional distributions $(Y|\phi)$ and $(\phi|Y)$ are easy to simulate, one can use a Gibbs sampler to simulate from the two-star model, by iteratively simulating from the conditional distributions till the Markov Chain converges. 
\end{remark}

 \subsection{Simulation Results}

\begin{figure*}[h]
\centering
\hbox{\hspace{3em}\begin{minipage}[c]{1.0\textwidth}
\includegraphics[width=1.8in]{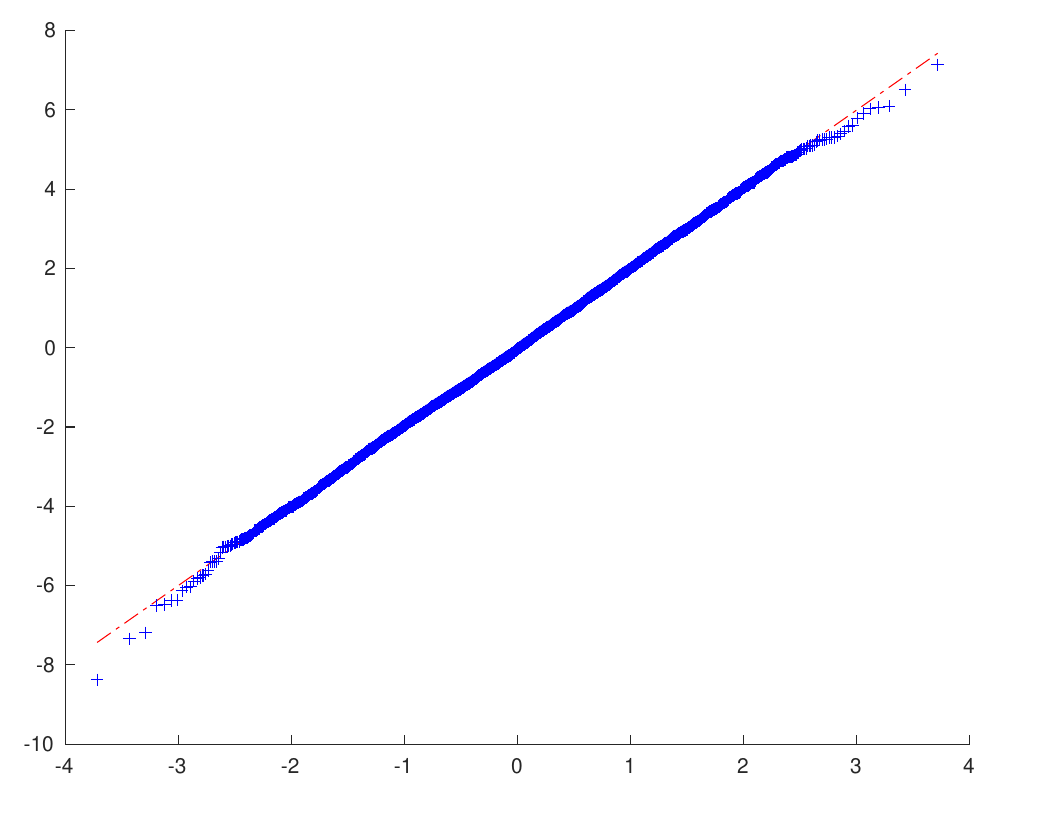}
\includegraphics[width=1.8in]{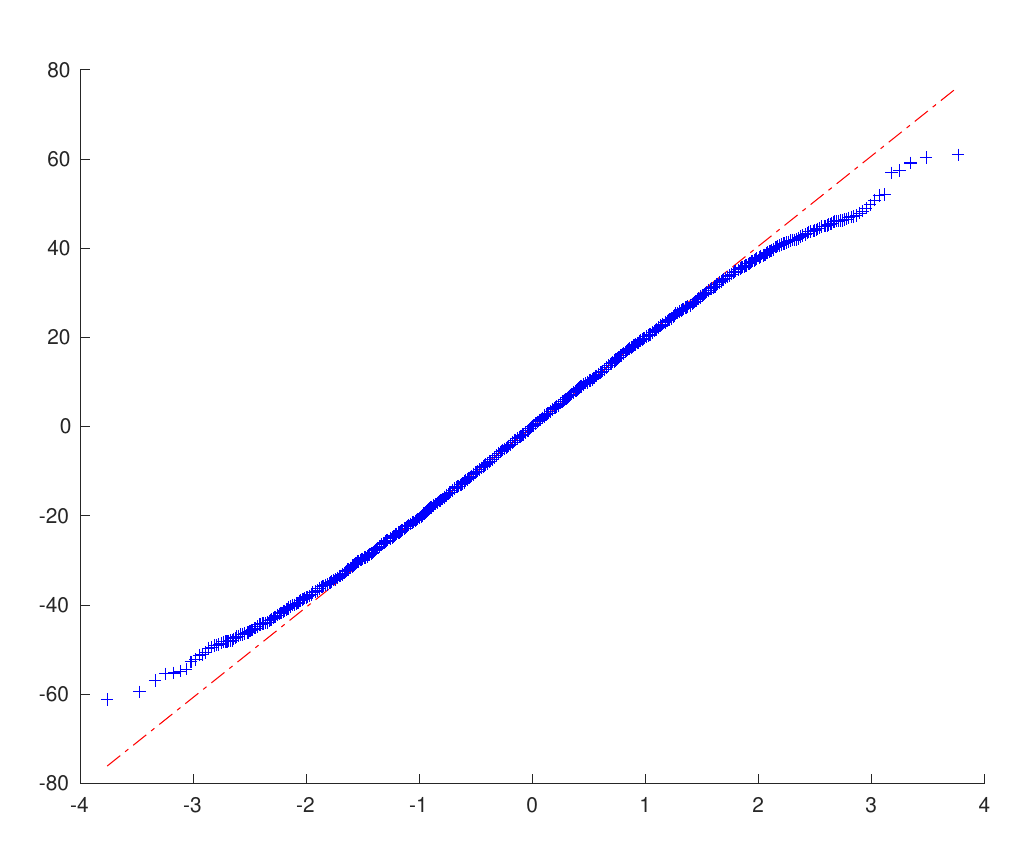}
 \includegraphics[width=1.8in]{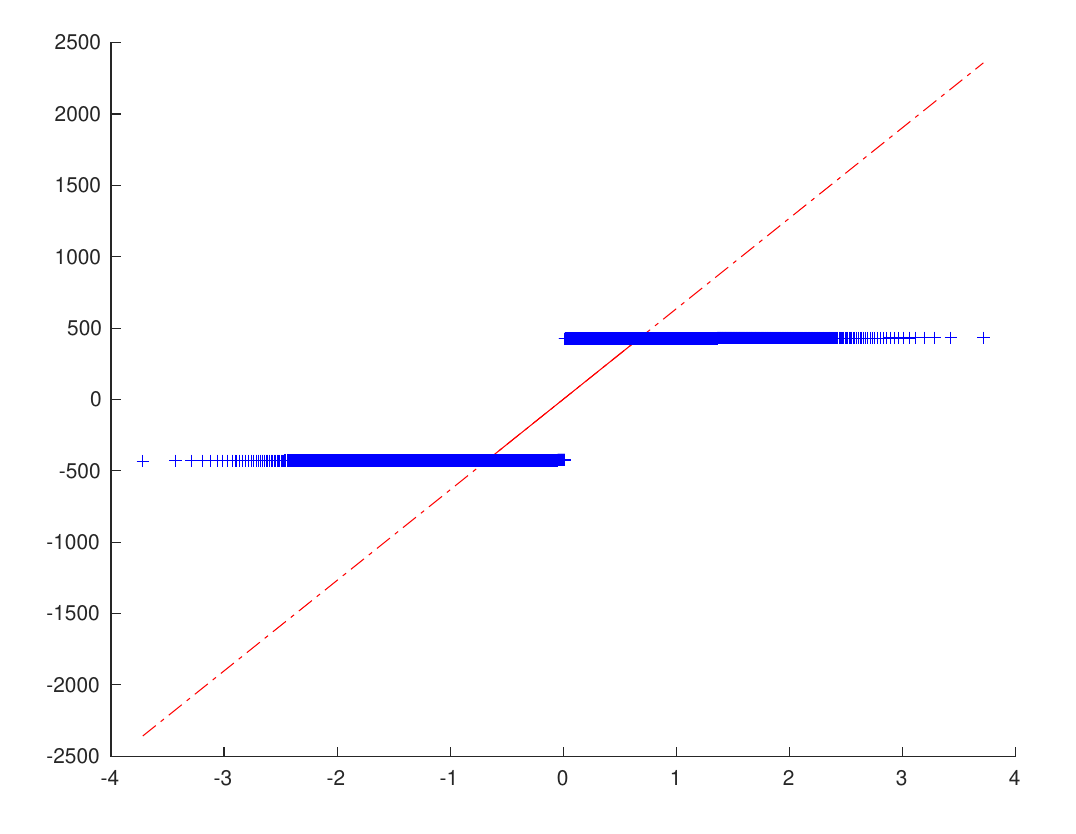}
\end{minipage}}
\caption{\footnotesize{The QQ plot for the centered and scaled sum of degrees is given for $5000$ independent samples from the two star ERGM on $n=500$ vertices, for the  three parameter configurations $(\theta,\beta)=(1/4,0)$, $(1/2,0)$, and $(3/4,0)$ respectively. In the uniqueness regime (first picture), the limiting distribution is Gaussian. At the critical point, the limiting distribution is no longer Gaussian. In the non uniqueness regime, the data is strongly bimodal. }}
\label{fig:fb1}
\end{figure*}

In this section  we validate Theorem \ref{mean1} and Theorem \ref{bm} using numerical simulations. For simulating from the two-star ERGM we use the Gibbs sampling algorithm of Proposition \ref{thm:bayes}. For verifying Theorem \ref{mean1}, we work with $n=500$ vertices on the two star ERGM with parameters $(\theta,\beta)$ equal to $(1/4, 0)$, $(1/2, 0)$, and  $(3/4,0)$, which belong to the uniqueness regime, the critical point, and the non-uniqueness regime respectively. For each of these three parameter configurations, we simulate $5000$ independent samples from the two star ERGM, by running the Gibbs sampling algorithm with a burn in period of $1000$ for each sample. For each sample, we observe the centered and scaled sum of degrees $$\frac{1}{n}\sum_{i=1}^nk_i(Y)=\frac{2}{n}\left[\sum_{i=1}^nd_i(G)-\frac{n(n-1)}{2}\right]=\frac{4}{n}\left[E(G)-\frac{n(n-1)}{4}\right].$$ The QQ plot of these values for the three regimes are given in figure \ref{fig:fb1}.

As is seen in figure \ref{fig:fb1}, in the uniqueness regime (first picture), the limiting distribution is clearly Gaussian, as there is a strong agreement with normal quantiles. At the critical point, the limiting distribution is no longer Gaussian, as is shown by deviation from the normal quantiles. In the non uniqueness regime, the data is strongly bimodal, and hence cannot be globally Gaussian. This is exactly the behavior predicted by Theorem \ref{mean1}. Theorem \ref{mean1} suggests that if we zoom into each of the two modes, we will again see Gaussian fluctuations. To confirm this, we do a QQ plot for the positive and negative values separately. This is given below in figure \ref{ref:fig2}.

\begin{figure*}[h]
\centering
\hbox{\hspace{5em}
\begin{minipage}[c]{1.0\textwidth}
\includegraphics[width=2in]{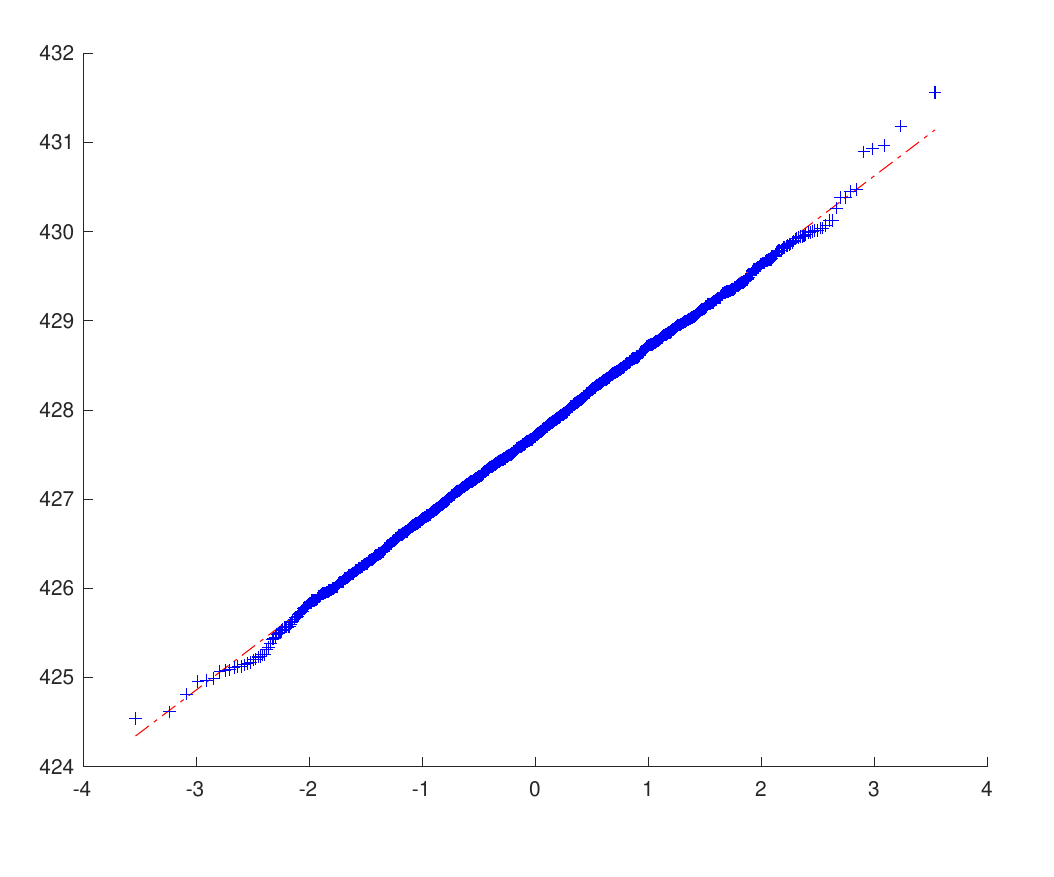}
\includegraphics[width=2in]{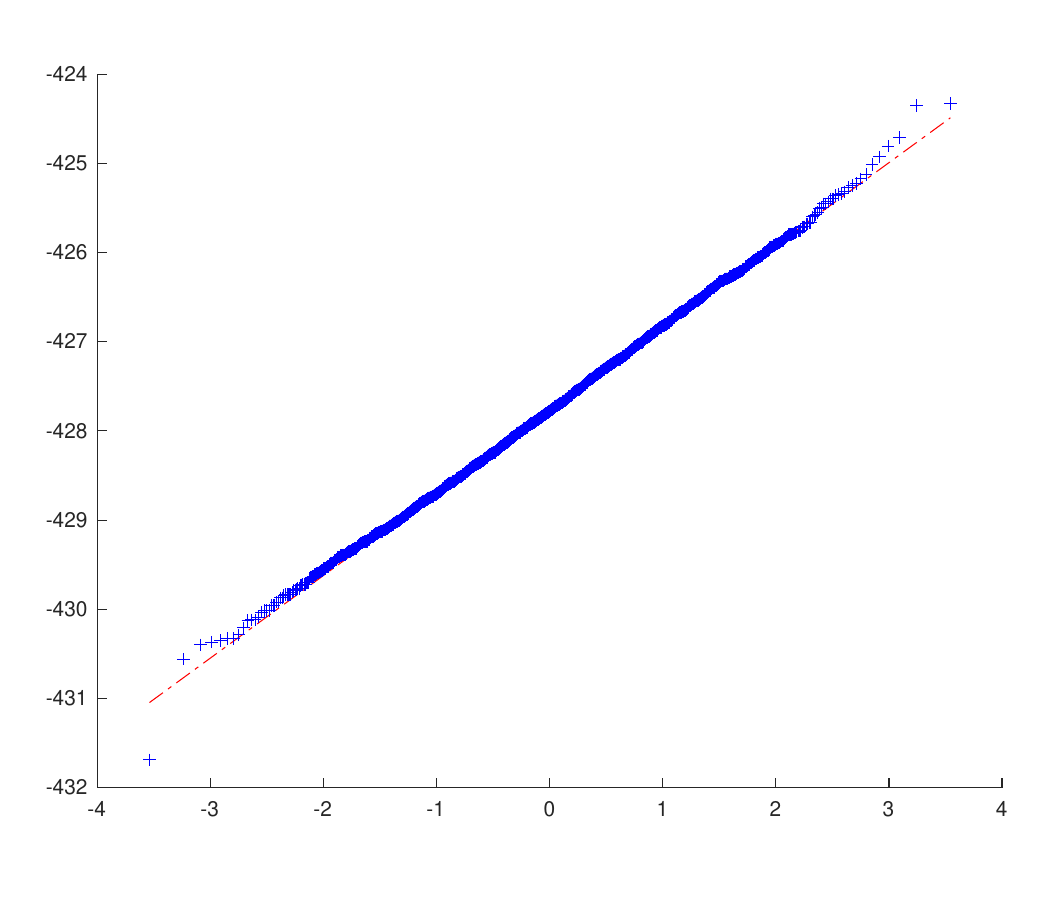}
\end{minipage}}
\caption{\footnotesize{Starting from $5000$ independent samples from the two star ERGM on $n=500$ vertices with parameter $(\theta,\beta)=(3/4,0)$, the QQ plot for the positive and negative values are given separately. Both the individual plots show agreement with Gaussian quantiles, which shows conditional Gaussian behavior near each of the two models. }}
\label{ref:fig2}
\end{figure*}

\begin{figure*}[h]
\centering
\includegraphics[width=5in,height=2in]{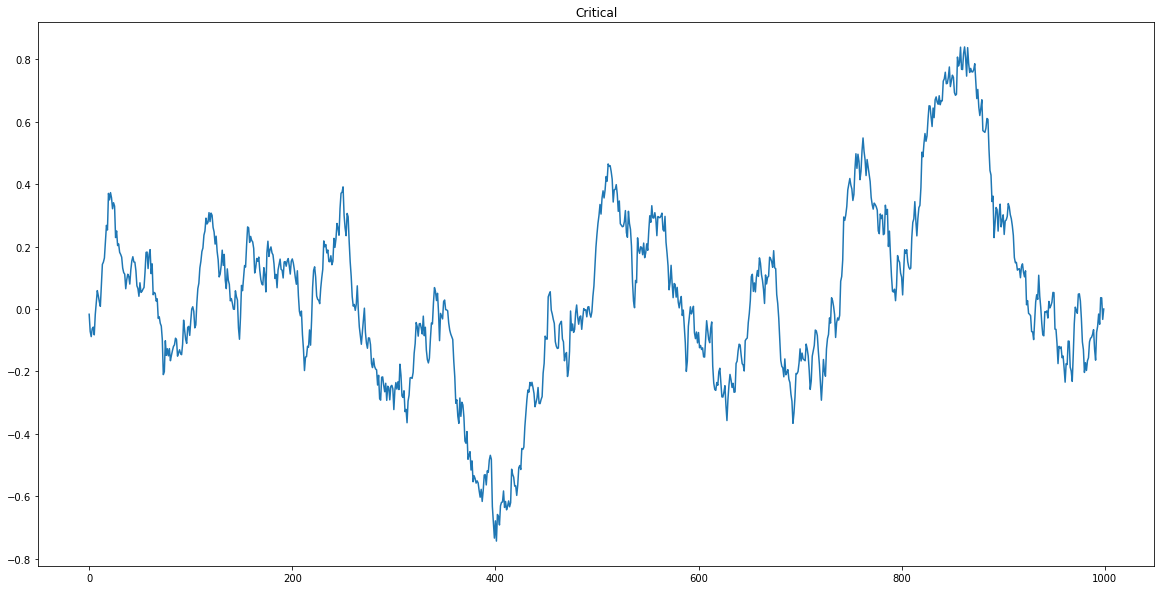}
 \caption{\footnotesize{From one sample from the two star ERGM on $n=1000$ vertices with parameter $(\theta,\beta)=(1/2,0)$, the centered  and partial sums of the degrees upto vertex $i$ is plotted against $i$, for $1\le i\le 1000$. The figure roughly resembles a Brownian curve starting and ending at the origin. }}
\label{ref:fig3}
\end{figure*}

For verifying Theorem \ref{bm}, we obtained one sample from the two star ERGM on $n=1000$ vertices at criticality ($(\theta,\beta)=(1/2,0)$), after running the chain for $1000$ iterations. Having obtained the graph $G$, we computed the partial sums 
$$\frac{1}{n-1}\sum_{j=1}^{i}\Big(k_j(Y)-\bar{k}(Y)\Big)=\frac{2}{n-1}\sum_{j=1}^i\Big(d_j(G)-\bar{d}(G)\Big),$$
and plotted the partial sums versus $i$ for $1\le i\le n$ in figure \ref{ref:fig3}.
 
As predicted, the plot looks like a Brownian curve starting and ending at $0$. Similar pictures were obtained in all parameter regimes.\\

The rest of the paper is as follows: Sections \ref{sec:two} proves Theorem \ref{mean1}, Theorem \ref{sd}, Corollary \ref{cor:est}, and Theorem \ref{bm}. The lemmas necessary for proving the main results are proved in section \ref{sec:three} for the uniqueness and non-uniqueness domains (i.e. $(\theta,\beta)\in \Theta_1\cup \Theta_2$), and in section \ref{sec:critical} for the critical domain (i.e. $(\theta,\beta)\in \Theta_3$). The appendix collects the proof of some helpful results.

\section{Proof of main results (Theorems \ref{mean1}, \ref{sd}, \ref{bm} and Corollary \ref{cor:est})}\label{sec:two}

For proving our main results we need the following lemmas, the proof of which is deferred to sections \ref{sec:three} and \ref{sec:critical} for $(\theta,\beta)\in \Theta_1\cup \Theta_2$ and $(\theta,\beta)\in \Theta_3$ respectively. 

 \begin{lem}\label{lem:uniqueness_aux}

 \begin{enumerate}
 \item[(a)]
 For $(\theta,\beta)\in \Theta_1$, we have
 
  $$n(\bar{\phi}-t)\stackrel{D}{\to}N\left(-\frac{2\theta t(1-t^2)}{[1-\theta(1-t^2)][1-2\theta(1-t^2)]},\frac{1}{\theta-2\theta^2(1-t^2)}\right).$$
 
 \item[(b)]
 For $(\theta,\beta)\in \Theta_2$, we have
  $$\Big[n(\bar{\phi}-t)\Big|\bar{\phi}>0\Big]\stackrel{D}{\to}N\left(-\frac{2\theta t(1-t^2)}{[1-\theta(1-t^2)][1-2\theta(1-t^2)]},\frac{1}{\theta-2\theta^2(1-t^2)}\right).$$
 
 \item[(c)]
 For $(\theta,\beta)\in \Theta_3$, we have
 $\sqrt{n}\bar{\phi}\stackrel{D}{\to}\zeta$,
 where $\zeta$ is a random variable on $\R$ with density proportional to $e^{-\frac{\zeta^2}{2}-\frac{\zeta^4}{24}}$.
 
 \end{enumerate}
 \end{lem}

\begin{lem}\label{lem:contrast}
 
Setting 
\begin{align}\label{eq:a1}
a_1:=\theta-\theta^2(1-t^2),
\end{align}
 for all $(\theta,\beta)\in \Theta$ we have the following conclusions:
\begin{enumerate}
\item[(a)]
$\sqrt{n}\Big[\sum\limits_{i=1}^n(\phi_i-\bar{\phi})^2-a_1^{-1}\Big]\stackrel{D}{\longrightarrow}N(0,2a_1^{-2}).
$

\item[(b)]
For any triangular array of real numbers $(c_{n}(i),c_{n}(2),\ldots, c_{n}(n))$ such that $\sum\limits_{i=1}^nc_{n}(i)=0$ and $n^{-1}\sum\limits_{i=1}^nc_{n}(i)^2\rightarrow1$, we have
$
\sum\limits_{i=1}^nc_{n}(i)\phi_i\stackrel{D}{\longrightarrow}N(0,a_1^{-1}).
$

\item[(c)]
Setting $S_i(\phi):=\sum_{j=1}^i(\phi_j-\bar{\phi})$, for every $\varepsilon>0$ we have 
$$\limsup_{\delta\to0}\limsup_{n\to\infty}\P_n(\max_{i,j\in [n]:|i-j|\le n\delta}|S_i(\phi)-S_j(\phi)|>\varepsilon)=0.$$

\end{enumerate}

\end{lem}

\subsection{Asymptotic notation}

Throughout the rest of the paper we will use the following notations. Let $\{r_n\}_{n\ge 1}$ and $\{s_n\}_{n\ge 1}$ be two sequences of positive real numbers. Then we will say
\begin{itemize}
\item
$r_n=o(s_n)$ if $\lim_{n\to\infty}\frac{r_n}{s_n}=0$,

\item
$r_n=O(s_n)$ or $r_n\lesssim s_n$ if $\limsup_{n\to\infty}\frac{r_n}{s_n}<\infty,$

\item
$r_n=\Omega(s_n)$ if $\liminf_{n\to\infty}\frac{r_n}{s_n}>0$.
\end{itemize}

If $\{R_n\}_{n\ge 1}$ and $\{S_n\}_{n\ge 1}$ are sequences of random variables, we will say
\begin{itemize}
\item
$R_n=o_P(S_n)$ if $\frac{R_n}{S_n}\stackrel{P}{\to}0$.

\item
$R_n=O_P(S_n)$, if $\frac{R_n}{S_n}$ is tight.
\end{itemize}

\subsection{ {Proof of Theorem \ref{mean1} }}

\begin{enumerate}
\item[(a)] $(\theta,\beta)\in \Theta_1$.

To begin, using \eqref{eq:phi_z} we have
 \[n(\bar{\phi}-t)=\frac{n}{n-1}\Big[\bar{k}(Y)-(n-1)t\Big]+\frac{n\bar{W}}{\sqrt{(n-1)\theta}}.\]
Using this along with part (a) of Lemma \ref{lem:uniqueness_aux} and the observation $\frac{n\bar{W}}{\sqrt{(n-1)\theta}}\stackrel{D}{\to}N(0, \frac{1}{\theta})$
gives
\[\bar{k}(Y)-(n-1)t\stackrel{D}{\longrightarrow}N\left(-\frac{2\theta t(1-t^2)}{[1-\theta(1-t^2)] [1-2\theta(1-t^2)]}, \frac{2(1-t^2)}{1-2\theta(1-t^2)}\right)\] Finally use \eqref{eq:K} to note that $\bar{k}(Y)=2\bar{d}(G)-(n-1)$, and so
\[\bar{d}(G)-(n-1)p\stackrel{D}{\longrightarrow}N\left(-\frac{\theta t(1-t^2)}{[1-\theta(1-t^2)] [1-2\theta(1-t^2)]}, \frac{1-t^2}{2[1-2\theta(1-t^2)]}\right),\]
which verifies Theorem \ref{mean1} for $(\theta,\beta)\in \Theta_1$.

\item[(b)]
By symmetry we have $\P_n(\bar{\phi}>0)=\P_n(\bar{\phi}<0)=\frac{1}{2}$, and on the set $\bar{\phi}>0$ $(\bar{\phi}<0)$ we have $\bar{\phi}\stackrel{P}{\to}t$ ($\bar{\phi}\stackrel{P}{\to}-t$) respectively, by invoking Lemma \ref{lem:uniqueness_aux} part (b).
 On the set $\bar{\phi}>0$, using \eqref{eq:phi_z} it follows that $$\frac{1}{n}\bar{k}(Y)\stackrel{P}{\to}t\Rightarrow \frac{2E(G)}{n^2}\stackrel{P}{\to} \frac{1+t}{2}=p>\frac{1}{2}.$$
 A similar argument gives that on the set $\bar{\phi}<0$ we have $$\frac{2E(G)}{n^2}\stackrel{P}{\to}\frac{1-t}{2}=1-p<\frac{1}{2}.$$
Thus $\P_n(\{\bar{\phi}>0\}\Delta\{\E(G)>\frac{n^2}{4}\})\to 0$ (here $\Delta$ represents symmetric difference between the two sets), and so without loss of generality we can replace the conditioning set $E(G)>\frac{n^2}{2}$ by $\bar{\phi}>0$. From then, using part (b) of Lemma \ref{lem:uniqueness_aux} and mimicking the proof of part (a) above gives the desired conclusion.

\item[(c)]

Again using \eqref{eq:phi_z} we have
$$\sqrt{n}\bar{\phi}=\frac{\sqrt{n}\bar{k}(Y)}{n-1}+\frac{\sqrt{n}\bar{W}}{\sqrt{(n-1)\theta}}.$$
Since $\bar{W}\stackrel{P}{\to}0$, it follows from part (b) of Lemma \ref{lem:uniqueness_aux} that $\frac{\bar{k}(Y)}{\sqrt{n}}\stackrel{D}{\to}\zeta$. The desired result then follows from on noting that
$$\frac{\bar{k}(Y)}{\sqrt{n}} =\sqrt{n}\Big[ \frac{2\bar{d}(G)-(n-1)}{n}\Big]=2\sqrt{n}\Big[\frac{\bar{d}(G)}{n^2}-\frac{1}{2}\Big]+O\Big(\frac{1}{\sqrt{n}}\Big).$$
\end{enumerate}

\subsection{ {Proof of Theorem \ref{sd} }}

Using \eqref{eq:phi_z} we can write 
\begin{align}\label{eq:phi_zsd}
\sum\limits_{i=1}^n(\phi_i-\bar{\phi})^2-a_1^{-1}
=A_n+B_n+C_n,
\end{align}
where 
\begin{align*}
A_n:=&\Big[\frac{1}{(n-1)\theta}\sum\limits_{i=1}^n(W_i-\bar{W})^2-\frac{1}{\theta}\Big],\quad
B_n:=\frac{2}{\sqrt{(n-1)^3\theta}}\sum\limits_{i=1}^n\Big(k_i(Y)-\bar{k}(Y)\Big)(W_i-\bar{W}),\\
C_n:=&\Big[\frac{1}{(n-1)^2}\sum\limits_{i=1}^n\Big(k_i(Y)-\bar{k}(Y)\Big)^2-\tau].
\end{align*}
Here we have used the fact that 
\begin{align}\label{eq:atau}
a_1^{-1}=\tau+\theta^{-1},
\end{align} 
where $a_1$ is as in \eqref{eq:a1}. We now claim that given the graph $G$, the random variables $\sqrt{n}A_n$ and $\sqrt{n}B_n$ are asymptotically independent, i.e. for any $s\in \R$ we have
\begin{align}\label{eq:ind_ab}
\Big|\E(e^{is\sqrt{n}(A_n+B_n)}|G)-\E(e^{is\sqrt{n}A_n}|G)\E(e^{is\sqrt{n}A_n}|G)\Big|\stackrel{P}{\to}0.
\end{align}
Given \eqref{eq:ind_ab}, we have
\begin{align}\label{useful?}
\mathbb{E}e^{is\sqrt{n}(A_n+B_n+C_n)}=\mathbb{E}[\mathbb{E}(e^{is\sqrt{n}A_n}|G) \mathbb{E}(e^{is\sqrt{n}B_n}|G)e^{is\sqrt{n}C_n}]+o_P(1).
\end{align}
Also, note that
$$A_n\stackrel{P}{\to}0,\quad {\sf Var}(B_n|G)=\frac{4}{(n-1)^3\theta}\sum_{i=1}^n\left(k_i(Y)-\bar{k}(Y)\right)^2\stackrel{P}{\to}0,$$
where the second convergence uses the observation $\sum_{i=1}^n(k_i-\bar{k})^2=O_P(n^2)$. This, along with \eqref{eq:phi_zsd} gives
$$\frac{1}{n^2}\sum_{i=1}^n\Big(k_i(Y)-\bar{k}(Y)\Big)^2\stackrel{P}{\to} \tau.$$
 Consequently, given $G$ the random variable $\sqrt{n}B_n$ has a Normal distribution with mean $0$, and variance $D_n$, where 
\begin{align}\label{eq:b|g}
D_n:=\frac{4n}{(n-1)^3\theta}\sum_{i=1}^n\Big(k_i(Y)-\bar{k}(Y)\Big)^2\stackrel{P}{\to}\frac{4\tau}{\theta}=:\sigma_2^2.
\end{align}
Finally, it is straightforward to check that   \begin{align}\label{eq:A}
\sqrt{n}A_n\stackrel{D}{\to}N(0,\sigma_1^2),\text{ where } \sigma_1^2:=2\theta^{-2}\Big.
\end{align}
Combining \eqref{useful?} along with \eqref{eq:b|g} and \eqref{eq:A} gives
$$\mathbb{E}e^{is\sqrt{n}(A_n+B_n+C_n)}=e^{-\frac{s^2}{2}(\sigma_1^2+\sigma_2^2)}\mathbb{E}[e^{is\sqrt{n}C_n}]+o_P(1).$$
Since $\sqrt{n}(A_n+B_n+C_n)\stackrel{D}{\to}N\Big(0,2a_1^{-2}\Big)$ by Lemma \ref{lem:contrast} part (a), it follows hat $\sqrt{n}C_n\stackrel{D}{\to}N(0,\sigma_3^2)$, where $$2a_1^{-2}=\sigma_1^2+\sigma_2^2+\sigma_3^2=2\theta^{-2}+4\tau \theta^{-1}+\sigma_3^2.$$ Using \eqref{eq:atau} it follows that $\sigma_3^2=2\tau^2$, as desired.

%
To complete the proof, it suffices to verify \eqref{eq:ind_ab}. To this effect, given the graph $G$ construct an orthogonal matrix ${\bf O}_n$ whose first row is proportional to the constant vector ${\bf 1}$, and second row is proportional to the contrast vector $(k_1-\bar{k},\ldots,k_n-\bar{k})$. Then with ${\bf U}:={\bf O}_n'{\bf W}\sim N({\bf 0},{\bf I}_n)$ we have $$\sum_{i=2}^nU_i^2=\sum_{i=1}^n(W_i-\bar{W})^2,\quad U_2=\frac{1}{\sqrt{\sum\limits_{i=1}^n(k_i-\bar{k})^2}}\sum\limits_{i=1}^n(k_i-\bar{k})W_i,$$
and so
\begin{align*}
\sqrt{n}(A_n,B_n)=&\bigg\{\sqrt{n}\Big[\frac{1}{(n-1)\theta}\sum\limits_{i=2}^nU_i^2-\frac{1}{\theta}\Big]
,U_2\sqrt{D_n}\bigg\}\\
=&\bigg\{\sqrt{n}\Big[\frac{1}{(n-1)\theta}\sum\limits_{i=3}^nU_i^2-\frac{1}{\theta}\Big]
,U_2\sqrt{D_n}\bigg\}+o_P(1),
\end{align*}
and so \eqref{eq:ind_ab} follows.

\subsection{Proof of Corollary \ref{cor:est}}
\begin{enumerate}
\item[(a)]
Let $$\hat{t}:=\frac{2 E(Y)}{n^2},\quad \hat{\tau}:=\frac{1}{n^2}\sum_{i=1}^n\Big(k_i(Y)-\bar{k}(Y)\Big)^2.$$
Then using Theorem \ref{mean1} and Theorem \ref{sd}, it follows that for all $(\theta,\beta)\in \Theta$ we have
$n(\hat{t}^2-t^2)=O_P(1)$, and $\sqrt{n}(\hat{\tau}-\tau)=O_P(1)$. 
Thus, with $\hat{\theta}:=\frac{1}{1-\hat{t}^2}-\frac{1}{\hat{\tau}}$ we have
\begin{align*}
|\hat{\theta}-\theta|\le &\Big|\frac{1}{1-\hat{t}^2}-\frac{1}{1-t^2}\Big|+\Big|\frac{1}{\tau}-\frac{1}{\hat{\tau}}\Big|\\
=&|\hat{t}^2-t^2|\frac{1}{(1-t^2)(1-\hat{t}^2)}+\frac{\tau-\hat{\tau}}{\tau\hat{\tau}}
=O_P(|\hat{t}^2-t^2|)+O_P(|\hat{\tau}-\tau|)=O_P(n^{-1/2}).
\end{align*}
Similarly, with $\hat{\beta}:=\arctanh(\hat{t})-2\hat{\theta}\hat{t}$ we have
\begin{align*}
|\hat{\beta}-\beta|\le& |\arctanh(\hat{t})-\arctanh(t)|+2\hat{\theta}|\hat{t}-t|+2|t|(|\hat{\theta}-\theta|)\\
= &O_P(|\hat{t}-t|)+O_P(|\hat{\theta}-\theta|)=O_P(n^{-1/2}).
\end{align*}

\item[(b)]
Let $\hat{\beta}:=\arctanh(\hat{t})-2\theta \hat{t}$. Then we have
\begin{align*}
|\hat{\beta}-\beta|\le |\arctanh(\hat{t})-\arctanh(t)|+2\theta|\hat{t}-t|=O_P(|\hat{t}-t|)=O_P(n^{-1}),
\end{align*}
as desired.

\end{enumerate}

\subsection{ {Proof of Theorem \ref{bm}}}
\begin{proof}
We first check the convergence of finite dimensional distributions. For the sake of simplicity we check it for $2$ dimensional distributions. Fixing $0<s_1<s_2<1$, it suffices to show that for any $(r_1,r_2)\in \R^2$,
\[r_1W_{n}(s_1)+r_2W_n(s_2)\stackrel{D}{\rightarrow}N(0,\tau\psi), \text{ where } \psi:=r_1^2s_1(1-s_1)+r_2^2s_2(1-s_2)+2r_1r_2s_1(1-s_2).\]
With  \[b_n(i):=\frac{(r_1+r_2)}{\sqrt{n}}1_{\{1\leq i<ns_1\}}+\frac{r_2}{\sqrt{n}}1_{\{ns_1<i\leq ns_2\}}\text{ and } c_n(i):=b_n(i)-\bar{b}_n,\] we have 
{\small{\begin{align}
\begin{split}\label{eq:fd}
r_1W_n(s_1)+r_2W_n(s_2)=&\frac{(r_1+r_2)}{n\sqrt{n}}\sum\limits_{1\leq i\leq ns_1}\Big(k_i(Y)-\bar{k}(Y)\Big)+\frac{r_2}{n\sqrt{n}}\sum\limits_{ns_1<i\leq ns_2}\Big(k_i(Y)-\bar{k}(Y)\Big)+O\left(\frac{1}{\sqrt{n}}\right)\\
=&\frac{1}{n}\sum\limits_{i=1}^nb_n(i)\Big(k_i(Y)-\bar{k}(Y)\Big)+O\left(\frac{1}{\sqrt{n}}\right)=\frac{1}{n}\sum\limits_{i=1}c_n(i)k_i+O\left(\frac{1}{\sqrt{n}}\right).
\end{split}
\end{align}} }
Since $\sum\limits_{i=1}^nc_n(i)=0$ and \[\sum\limits_{i=1}^nc_n(i)^2=\sum\limits_{i=1}^nb_n(i)^2-n\bar{b}_n^2\rightarrow (r_1+r_2)^2s_1+r_2^2(s_2-s_1)-(r_1s_1+r_2s_2)^2=\psi,\]
by part (c) of Lemma \ref{lem:contrast} we have $\sqrt{n}\sum\limits_{i=1}^nc_n(i)\phi_{i}\stackrel{D}{\rightarrow}N\Big(0,\frac{\psi}{a_1}\Big).$ This, along with \eqref{eq:phi_z}, gives  $\frac{1}{n}\sum\limits_{i=1}^nc_n(i)k_i\stackrel{D}{\rightarrow}N(0,\tau\psi)$. This, along with \eqref{eq:fd} verifies convergence of finite dimensional distributions.

 It thus suffices to show tightness, for which using Arzela-Ascoli Theorem it suffices to verify that for every $\varepsilon>0$ we have
 \begin{align}\label{eq:tight}
 \limsup_{\delta\to0}\limsup_{n\to\infty}\P_n\Big(\max_{i,j\in [n]:|i-j|\le n\delta }|S_i({\bf d})-S_j({\bf d})|>(n-1)\varepsilon \Big)=0.
 \end{align}
 To verify \eqref{eq:tight}, first use \eqref{eq:phi_z} to note that
\begin{align*}
&\frac{1}{n-1}\max_{i,j\in [n]:|i-j|\le n\delta}|S_i({\bf d})-S_j({\bf d})|\\
\le& \frac{1}{2}\left[\max_{i,j\in [n]:|i-j|\le n\delta}|S_i(\phi)-S_j(\phi)|+\frac{1}{\sqrt{(n-1)\theta}}\max_{i,j\in [n]:|i-j|\le n\delta} |S_i({\bf W})-S_j({\bf W})|\right],
\end{align*}
 where ${\bf W}=(W_1,\ldots,W_n)$ is a sequence of i.i.d.~$N(0,1)$ random variables, and $S({\bf W})=\sum_{j=1}^iW_j$. This in turn gives the following bound to the RHS of \eqref{eq:tight}:
  \begin{align*}
&  \P_n\Big(\max_{i,j\in [n]:|i-j|\le n\delta}|S_i({\bf d})-S_j({\bf d})|>\varepsilon (n-1)\Big)\\
\le &\P_n\left( \max_{i,j\in [n]:|i-j|\le n\delta}|S_i(\phi)-S_j(\phi)|>\frac{\varepsilon}{4}\right)\\
+&\P_n\left(\frac{1}{\sqrt{(n-1)\theta}}\max_{i,j\in [n]:|i-j|\le n\delta} |S_i({\bf W})-S_j({\bf W})|>\frac{\varepsilon\sqrt{(n-1)\theta}}{4}\right).
\end{align*}

 The first term in the RHS above converges to $0$ as $n\to\infty$ followed by $\delta\to 0$ using part (c) of Lemma \ref{lem:contrast}, and the second term converges to $0$ under the same double limit by tightness of sample paths for partial sums of i.i.d.~random variables. Thus we have verified \eqref{eq:tight}, and hence the proof of the theorem is complete.
%
%
%
%
%
%
\end{proof}

\section{Proof of Lemma \ref{lem:uniqueness_aux} and Lemma \ref{lem:contrast} for $(\theta,\beta)\in \Theta_1\cup \Theta_2$\label{sec:three}}

We first state a general  approximation result, which will be used to analyze the marginal distribution of $(\phi_1,\ldots,\phi_n)$ by approximating the un-normalized density $f_n(.)$ of Proposition \eqref{thm:bayes} by something more tractable. The approximating measure will change across the three parameter regimes $\Theta_1\cup\Theta_2\cup\Theta_3$.

\begin{lem}\label{ctv}
For an interval $U\subseteq \R$, let $h_n(.),g_n(.):U^n\mapsto \R$ be non negative and integrable. Define the probability measures $\mathbb{G}_n$ and $\mathbb{H}_n$ on $U^n$ by setting
$$\frac{d\mathbb{G}_n}{d\lambda^{\otimes n}}:=\frac{g_n}{\int_{U^n}g_nd\lambda^{\otimes n}},\quad \frac{d\mathbb{H}_n}{d\lambda^{\otimes n}}:=\frac{h_n}{\int_{U^n}h_nd\lambda^{\otimes n}},$$
where $\lambda^{\otimes n}$ is Lebesgue measure on $\R^n$. 
Setting $L_n(.)=\log\frac{g_n}{h_n}$, suppose that 
$L_n$ is $O_P(1)$ under both measures $\G_n,\mathbb{H}_n$.

\begin{enumerate}
\item[(a)]
Then the sequence of probability measures $\G_n$ and $\mathbb{H}_n$ are mutually contiguous. 

\item[(b)]
If $(X_n,L_n)\stackrel{D,\mathbb{G}_n}{\longrightarrow}N(\mu_1,\mu_2,\sigma_1^2,\sigma_2^2,\sigma_{12})$ then $\mu_2+\frac{1}{2}\sigma_2^2=0$, and  $X_n\stackrel{D,\mathbb{H}_n}{\longrightarrow}N(\mu_1+\sigma_{12},\sigma_1^2).$

\item[(c)]
If  $L_n\stackrel{D,\mathbb{G}_n}{\rightarrow}c$ where $c$ is a constant, then $\lVert \mathbb{G}_n-\mathbb{H}_n\rVert_{TV}\rightarrow 0.$
\end{enumerate}
\end{lem}

 Our plan is use Lemma \ref{ctv} to approximate the distribution of $\phi$ by a multivariate Gaussian distribution. The following lemma summarizes some estimates under the approximating Gaussian distribution.

\begin{lem}\label{lem:G1}
Let $\mathbb{G}_{1n}$ be a multivariate Gaussian distribution on $\R^n$, with density proportional to $g_{1n}$, where  $$-\log g_{1n}(\phi)=\frac{n(n-1)}{2}p(t,t)+\frac{a_1n}{2}\sum_{i=1}^n(\phi_i-t)^2-\frac{a_2n^2}{2}(\bar{\phi}-t)^2,$$
with $a_1=\theta-\theta^2(1-t^2)$ as in \eqref{eq:a1}, and $a_2:=\theta^2(1-t^2)$.
Then the following conclusions hold under $\G_{1n}$.

\noindent(a)
$
\E_{\mathbb{G}_{1n}}|\phi_i-t|^\ell\lesssim_\ell n^{-\ell/2}.
$
\\

\noindent(b)
$
\sum\limits_{1\le i<j\le n}\Big(\phi_i+\phi_j-2t\Big)^4\stackrel{P}{\to} \frac{6}{a_1^2},
$
\\

\noindent(c)
Suppose ${\bf c}=(c_n(1),\ldots,c_n(n))$ be a vector such that $\sum_{i=1}^nc_n(i)=0$, and $\frac{1}{n}\sum_{i=1}^nc_n(i)^2\to 1$. Then we have 
{\small{\begin{equation*}
\left[n(\bar{\phi}-t),\sqrt{n}\Big(\sum\limits_{i=1}^{n}(\phi_{i}-\bar{\phi})^{2}-a_{1}^{-1}\Big),n\sum\limits_{i=1}^{n}(\phi_{i}-t)^{3},\sum_{i=1}^nc_i\phi_i\right]\stackrel{d}{\rightarrow}N({\bf 0},\Sigma)
\end{equation*}}
where   
\begin{equation*}\Sigma:=  \left[ \begin{array}{cccc}
   \frac{1}{a_1-a_2}& 0 & \frac{3}{a_1(a_1-a_2)}&0 \\
   0& \frac{2}{a_1^2} & 0&0 \\
   \frac{3}{a_1(a_1-a_2)}&0& \frac{15a_{1}-6a_{2}}{a_{1}^{3}(a_{1}-a_{2})}&0\\
   0&0&0&\frac{1}{a_1}
 \end{array}  \right].
\end{equation*}}

\noindent(d)
For every $\varepsilon>0$, setting $S_i(\phi)=\sum_{j=1}^i(\phi_j-\bar{\phi})$ as before, we have
$$\limsup_{\delta\to0}\limsup_{n\to\infty}\P_n(\max_{i,j\in [n]:|i-j|\le n\delta}|S_i(\phi)-S_j(\phi)|>\varepsilon )=0.$$

\end{lem}

The final result we need for proving Lemma \ref{lem:uniqueness_aux} in the regime $(\theta,\beta)\in \Theta_1\cup \Theta_2$ is the following:
\begin{lem}\label{lem:12whole}
Let $U=\R$ if $(\theta,\beta)\in \Theta_1$, and $U=(0,\infty)$ if $(\theta,\beta)\in \Theta_2$.

\begin{enumerate}
\item[(a)]
Then exists positive constants  $\lambda_1\ge \lambda_2$ such that for all $(x,y)\in U^2$ we have
\begin{align*}
\frac{\lambda_2}{2}[(x-t)^2+(y-t)^2]\leq p(x,y)-p(t,t)\leq\frac{\lambda_1}{2}[(x-t)^2+(y-t)^2] ,
\end{align*}
 
\item[(b)]
There exists $M$ large enough such that
\begin{align*}
\log \mathbb{P}_{n,U}(\sum\limits_{i=1}^n(\phi_i-t)^2> M)\lesssim -n.
\end{align*}
where $\mathbb{P}_{n,U}$ denotes the conditional law of $\phi$ under $\mathbb{P}_{n}$ given $\phi\in U^n$. 


\item[(c)]
For any $l\in\mathbb{N}$ we have
$
\mathbb{E}_{n,U}|\phi_i-t|^l\lesssim_l n^{-\ell/2},
$

\item[(d)]
$\E_{{n,U}}\big[\sum\limits_{i=1}^n(\phi_i-t)\big]^{2}\lesssim1$.
\end{enumerate}
\end{lem}

The proof of the three Lemmas \ref{ctv}, \ref{lem:G1} and \ref{lem:12whole} are deferred to the the appendix (section \ref{sec:four2}).

\subsection{Proof of Lemma \ref{lem:uniqueness_aux} and Lemma \ref{lem:contrast} for $(\theta,\beta)\in \Theta_1$}

Let $\mathbb{F}_n$ denote the marginal distribution of $\phi$ on $\R^n$ under $\P_n$, i.e. $\mathbb{F}_n$ is induced by the unnormalized density $f_n(.)$ defined in Proposition \ref{thm:bayes}. We begin by showing the following proposition:

\begin{ppn}\label{ppn:cont1}
The probability measures $\F_n$ and $\G_{1n}$ are mutually continuous, where $\G_{1n}$ is the multivariate Gaussian distribution introduced in Lemma \ref{lem:G1}. 
\end{ppn}

\begin{proof}
To this effect, with $q(.)$ as in Lemma \ref{lem:min}, use a Taylor's series expansion to get
 $$q\Big(\frac{x+y}{2}\Big)=q(t)+\frac{q''(t)}{2}\Big(\frac{x+y}{2}-t\Big)^2+\frac{q'''(t)}{3!}\Big(\frac{x+y}{2}-t\Big)^3+\frac{q''''(t)}{4!}\Big(\frac{x+y}{2}-t\Big)^4+R(x,y),$$
 where $|R(x,y)|\lesssim |x-t|^5+|y-t|^5$. Recalling that $p(x,y)=q\Big(\frac{x+y}{2}\Big)+\frac{\theta}{4}(x-y)^2$ then gives
 \begin{align*}
 p(x,y)=&\frac{\theta}{4}(x-y)^2+q(t)+\frac{q''(t)}{2}\Big(\frac{x+y}{2}-t\Big)^2+\frac{q'''(t)}{3!}\Big(\frac{x+y}{2}-t\Big)^3\\
 +&\frac{q''''(t)}{4!}\Big(\frac{x+y}{2}-t\Big)^4+R(x,y)\\
 =&p(t,t)+\frac{1}{2}\Big[a_1(x-t)^2+a_1(y-t)^2-2a_2(x-t)(y-t)\Big]+\frac{a_3}{3!}({x+y}-2t)^3\\
  +&\frac{a_4}{4!}({x+y}-2m)^4+R(x,y),
  \end{align*}
where $a_1=\theta-\theta^2(1-t^2), a_2=\theta^2(1-t^2)$ as in Lemma \ref{lem:G1} , and $ a_3:=\frac{q'''(t)}{8}, a_4:=\frac{q''''(t)}{16}.$
%
Adding $p(\phi_i,\phi_j)$ over $i<j$, this gives 
$$-\log f_n(\phi)=\sum_{\ell=1}^4 R_{\ell,f_n}+\sum_{i<j}R(\phi_i,\phi_j),$$
where
\begin{align}
\begin{split}\label{eq:f_r0}
R_{1,f_n}:= &  \sum\limits_{i<j}\frac{a_1}{2}[(\phi_i-t)^2+(\phi_j-t)^2]=\frac{a_1(n-1)}{2}\sum\limits_{i=1}^n(\phi_i-t)^2 \\ 
  R_{2,f_n}:=&-a_2\sum\limits_{i<j}(\phi_i-t)(\phi_j-t)  =  -\frac{a_2n^2}{2}(\bar{\phi}-t)^2+\frac{a_2}{2}\sum\limits_{i=1}^n(\phi_i-t)^2\\
  R_{3,f_n}:=&\frac{a_3}{6}\sum\limits_{i<j}(\phi_i+\phi_j-2t)^3=\frac{a_3}{12}\left[\sum\limits_{i,j=1}^n(\phi_i+\phi_j-2t)^3-8\sum\limits_{i=1}^n(\phi_i-t)^3\right]\\
  =&\frac{(n-4)a_3}{6}\sum\limits_{i=1}^n(\phi_i-t)^3+\frac{3na_3}{6}(\bar{\phi}-t)\sum\limits_{i=1}^n(\phi_i-t)^2,\\
  R_{4,f_n}:=&\frac{a_4}{4!}\sum\limits_{1\leq i<j\leq n}(\phi_i+\phi_j-2t)^4.
  \end{split}
  \end{align}
  Consequently we have
  \begin{align}
  \begin{split}\label{eq:f_r}
  &\left|\log f_n(\phi)-\frac{a_1n}{2}\sum_{i=1}^n(\phi_i-t)^2-\frac{a_3n}{6}\sum_{i=1}^n(\phi_i-t)^3\right|\\
\lesssim
&n^2(\bar{\phi}-t)^2+n(\bar{\phi}-t)\sum_{i=1}^n(\phi_i-t)^2+\sum_{\ell=2}^5\sum_{i=1}^n|\phi_i-t|^\ell
\end{split}
\end{align}
      Fixing $b_4>a_3^2/3a_1$, define the function $h_{1n}(\phi)$ by
   \begin{align}
   \begin{split}\label{eq:h_r}
   -\log h_{1n}(\phi):=& \frac{n(n-1)}{2}p(t,t)+ \frac{a_{1}n}{2}\sum\limits_{i=1}(\phi_i-t)^2+ \frac{a_{3}n}{3!}\sum\limits_{i=1}^n(\phi_i-t)^3+\frac{b_{4}n}{4!}\sum\limits_{i=1}^n(\phi_i-t)^4\\
   =&\frac{n(n-1)}{2}p(t,t) +n^2 \eta(\phi_i-t),\quad \eta(x):=\frac{a_1}{2!}x^2+\frac{a_3}{3!}x^4+\frac{b_4}{4!}x^4,
   \end{split}
   \end{align}
   and note that $(\phi_1-t,\ldots,\phi_n-t)$ are i.i.d.~under $\mathbb{H}_{1n}$ with density proportional to $e^{-n^2\eta(.)}$, where $\mathbb{H}_{1n}$ denotes the probability measure induced by $h_{1n}$.
   It follows from straightforward calculus that
   \begin{align*}
\Big|\int_{\R} x^\ell e^{-x^2\eta(x)}dx\Big|\lesssim_\ell \frac{1}{n^{\frac{\ell+1}{2}}}\text{ if }\ell\text{ is even},\\
\lesssim_\ell \frac{1}{n^{\frac{\ell+3}{2}}} \text{ if }\ell\text{ is odd},
\end{align*}
and so
\begin{align}
\begin{split}\label{eq:prop2}
\E_{\mathbb{H}_{1n}}(\phi_i-t)^\ell \lesssim &\frac{1}{n^{\frac{\ell}{2}}}\text{ if }\ell\text{ is even},\\
\lesssim_\ell &\frac{1}{n^{\frac{\ell}{2}+1}} \text{ if }\ell\text{ is odd}.
\end{split}
\end{align}
 %
  Also, comparing \eqref{eq:f_r} and \eqref{eq:h_r} we have
  \begin{align}\label{eq:fh_r}
  |\log f_n(\phi)-\log h_{1n}(\phi)|\lesssim n^2(\bar{\phi}-t)^2+\Big|n(\bar{\phi}-t)\sum_{i=1}^n(\phi_i-t)^2\Big|+\sum_{\ell=2}^5\sum_{i=1}^n|\phi_i-t|^\ell.
  \end{align}
Using parts (c) and (d) of Lemma \ref{lem:12whole} it follows that $\log f_n(\phi)-\log h_{1n}(\phi)$ is $O_P(1)$ under $F_{n}$. To show the same conclusion under $\mathbb{H}_{1n}$, it suffices to note that 
\begin{align}\label{eq:h_est1}
 \E_{\mathbb{H}_{1n}}\Big[\sum_{i=1}^n(\phi_i-t)\Big]^2\lesssim 1,\quad \E|\phi_i-t|^\ell \lesssim n^{-\ell/2},
 \end{align}
both of which follow from
\eqref{eq:prop2}. It thus follows from Lemma \ref{ctv} that $\F_n$ and $\mathbb{H}_{1n}$ are mutually contiguous. To complete the proof, it suffices to show that $\G_{1n}$ and $\mathbb{H}_{1n}$ are mutually contiguous. Proceeding to verify this, note that
\begin{align}
\begin{split}\label{eq:gh_r}
\left|\log \frac{g_{1n}(\phi)}{h_{1n}(\phi)}\right|\lesssim  &n(\bar{\phi}-t)^2+ \Big|n(\bar{\phi}-t)\sum_{i=1}^n(\phi_i-t)^2\Big|\\
+&n\Big|\sum_{i=1}^n(\phi_i-t)^3\Big|+\sum_{\ell=2}^5\sum_{i=1}^n|\phi_i-t|^\ell.
\end{split}
\end{align}
We need to show that the RHS of \eqref{eq:gh_r} is $O_P(1)$ under both $\mathbb{H}_{1n}$ and $\G_{1n}$. Again the desired conclusion for $\mathbb{H}_{1n}$ follows \eqref{eq:h_est1}, and using \eqref{eq:prop2} to note that
\begin{align}\label{eq:h_est2}
n\E_{\mathbb{H}_{1n}}\Big[\sum_{i=1}^n(\phi_i-t)^3\Big]^2\lesssim 1.
\end{align}
To complete the proof, it suffices to verify \eqref{eq:h_est1} and \eqref{eq:h_est2} under $\G_{1n}$. But this follows from parts (a) and (c) of Lemma \ref{lem:G1}. This shows that $\F_n$ and $\G_{1n}$ are mutually continuous, and so we have verified the proposition.
\end{proof}

\begin{proof}[Proof of Lemma \ref{lem:uniqueness_aux} for $(\theta,\beta)\in \Theta_1$]
Use \eqref{eq:f_r0} to note that
\begin{align*}
-\log\frac{f_n}{g_{1n}}=\frac{a_2-a_1}{2}\sum_{i=1}^n(\phi_i-t)^2+R_{3,f_n}+R_{4,f_n}+\sum_{1\le i<j\le n}R(\phi_i,\phi_j).
\end{align*}
Invoking parts (a) and (b) of Lemma \ref{lem:G1}, under $\G_{1n}$ we have
\begin{align}\label{eq:three}
\left(\sum\limits_{i=1}^n(\phi_i-t)^2,R_{4,f_n},\sum\limits_{1\leq i<j\leq n}R(\phi_i,\phi_j)\right)\stackrel{P}{\longrightarrow}\Big(\frac{1}{a_1},\frac{a_4}{4a_1^2},0\Big).
\end{align}
Also, using \eqref{eq:f_r0}, a direct expansion gives
\begin{align}\label{eq:r3}
\notag R_{3,f_n}=&\frac{(n-4)a_3}{6}\sum_{i=1}^n(\phi_i-t)^3+\frac{3a_3}{6}n(\bar{\phi}-t)\sum_{i=1}^n(\phi_i-t)^2\\
=&\frac{na_3}{6}\sum_{i=1}^n(\phi_i-t)^3+\frac{a_3}{2a_1}n(\bar{\phi}-t)+o_p(1),
\end{align}
where the last equality again uses \eqref{eq:three}. Combining \eqref{eq:three} and \eqref{eq:r3} along with \eqref{eq:f_r0} gives
\begin{align}\label{eq:log_like}
-\log \frac{f_n}{g_{1n}}=\frac{a_2-a_1}{2a_1}+\frac{a_4}{4a_1^2}+\frac{na_3}{6}\sum_{i=1}^n(\phi_i-t)^3+\frac{a_3}{2a_1}n(\bar{\phi}-t)+o_p(1)
\end{align}
Now using part (c) of Lemma \ref{lem:G1}, under $\G_{1n}$ we have \[\left[n(\bar{\phi}-t), n\sum_{i=1}(\phi_i-t)^3\right]\stackrel{d}{\to}N\left({\bf 0},  \left[ \begin{array}{ccc}
   \frac{1}{a_1-a_2}&  \frac{3}{a_1(a_1-a_2)} \\
   \frac{3}{a_1(a_1-a_2)}&\frac{15a_{1}-6a_{2}}{a_{1}^{3}(a_{1}-a_{2})}
 \end{array}  \right]\right).\]
Using the fact that $\F_n$ and $\G_{1n}$ are contiguous, it follows from part (b) of Lemma \ref{ctv} that under  $\F_n$ we have $n(\bar{\phi}-t)\stackrel{d}{\to}N(-\mu, \frac{1}{a_1-a_2})$, where
\[\mu=\frac{a_3}{2a_1}\times \frac{1}{a_1-a_2}+\frac{a_3}{6}\times \frac{3}{a_1(a_1-a_2)}=\frac{a_3}{a_1(a_1-a_2)}=\frac{2\theta t(1-t^2)}{[1-\theta(1-t^2)] [1-2\theta(1-t^2)]}\]
as desired.

\end{proof}

\begin{proof}[Proof of Lemma \ref{lem:contrast} for $(\theta,\beta)\in \Theta_1$]
\begin{enumerate}
\item[(a)]
Using part (c) of Lemma \ref{lem:G1} along with \eqref{eq:log_like} it follows that the random variables $\sqrt{n}[\sum_{i=1}^n(\phi_i-t)^2-\frac{1}{a_1}]$ and $\frac{\log f_n}{\log g_{1n}}$ are asymptotically mutually independent and Gaussian under $\G_{1n}$. The desired result then follows from part (c) of Lemma \ref{lem:G1} along with part (b) of Lemma \ref{ctv}.


\item[(b)]
The proof of part (b) follows on similar lines as the proof of part (a), and is not repeated here.

\item[(c)]
By Proposition \ref{ppn:cont1} the two distributions $\F_n$ and $\G_{1n}$ are mutually contiguous, and so it suffices to verify the result under $\G_{1n}$. But this is precisely part (d) of Lemma \ref{lem:G1}, and so the proof is complete.

\end{enumerate}
\end{proof}

\subsection{Proof of Lemma \ref{lem:uniqueness_aux} and Lemma \ref{lem:contrast} for $(\theta,\beta)\in \Theta_2$} 

We begin by stating the following proposition, the proof of which is deferred to the appendix \ref{sec:four3}.

\begin{ppn}\label{lem:sep}
For $(\theta,\beta)\in \Theta_2$, we have $$\left|\frac{1}{2}-\P_n(\phi_i\ge 0,1\le i\le n)\right|\le e^{-\Omega(n)}.$$
\end{ppn}

\begin{proof}[Proof of Lemma \ref{lem:uniqueness_aux} and Lemma \ref{lem:contrast} for $(\theta,\beta)\in \Theta_2$]

By symmetry, we have $\mathbb{P}_n(\bar{\phi}>0)=\frac{1}{2}$, which along with Proposition \ref{lem:sep} gives that conditioned on $\bar{\phi}>0$ we have 
$$\mathbb{P}_n(\phi_i\leq 0 \text{ for some }i, 1\leq i\leq n|\bar{\phi}>0)\leq e^{-\Omega(n)}.$$
Thus at an exponentially vanishing cost we can replace the event $\bar{\phi}>0$ by  the event $\{\phi_i>0,1\leq i\leq n\}$. Consequently, invoking Lemma \ref{lem:12whole} with  $U=(0,\infty)$ and proceeding exactly as in the uniqueness domain we get the following conclusions:
\begin{align*}
&[n(\bar{\phi}-t)|\bar{\phi}>0]\stackrel{D}{\rightarrow}N\left(-\frac{2\theta t(1-t^2)}{[1-\theta(1-t^2)][1-2\theta(1-t^2)]},\frac{1}{\theta-2\theta^2(1-t^2)}\right),\\
&\left(\sqrt{n}\Big[\sum\limits_{i=1}^n(\phi_i-\bar{\phi})^2-a_1^{-1}\Big]\Big|\bar{\phi}>0\right)\stackrel{D}{\rightarrow}N(0,2a_1^{-2}),\\
&\left(\sqrt{n}\sum\limits_{i=1}^nc_n(i)(\phi_i-\bar{\phi})|\bar{\phi}>0\right)\stackrel{D}{\rightarrow}N(0,a_1^{-1}),\\
&\limsup_{\delta\to 0}\limsup_{n\to\infty}\P_n(\max_{i,j\in [n]:|i-j|\le n\delta}|S_i(\phi)-S_j(\phi)|>\varepsilon)=0.
\end{align*}
Here the last line above holds for any $\varepsilon>0$. Similarly calculations hold on the set $\bar{\phi}<0$ as well.
This readily proves Lemma \ref{lem:uniqueness_aux}. Lemma \ref{lem:contrast} follows on noting that the conditional distribution in the second, third and fourth lines in the above display is the same for $\bar{\phi}>0$ and $\bar{\phi}<0$.
\end{proof}

\section{Proof of Lemmas \ref{lem:uniqueness_aux} and \ref{lem:contrast} for $(\theta,\beta)\in \Theta_3$}\label{sec:critical}




We first state two lemmas which we will use to prove Lemma \ref{lem:uniqueness_aux} and Lemma \ref{lem:contrast} for $(\theta,\beta)\in \Theta_3$. The first lemma is the analogue of Lemma \ref{lem:12whole} parts (c) and (d), and the second lemma is the analogue of Lemma \ref{lem:G1}.  The proof of the two lemmas are deferred to the appendix (\ref{sec:four4}).
\begin{lem}\label{l31}
Suppose $(\theta,\beta)\in \Theta_3$.
\begin{enumerate}
\item[(a)]
For any positive integer $\ell$ we have
\begin{align*}
\mathbb{E}|\phi_i-\bar{\phi}|^l\lesssim_\ell \frac{1}{n^{l/2}}.
\end{align*}
%
%
%
%
\item[(b)]
$\limsup_{n\to\infty}n^2\E \bar{\phi}^4<\infty.$
\end{enumerate}
\end{lem}



\begin{lem}\label{lem:G2}
Suppose $$g_{3n}(\phi):=\frac{(n-1)\theta}{4}\sum\limits_{i=1}^n(\phi_i-\bar{\phi})^2-\frac{1}{2}n\bar{\phi}^2-\frac{1}{24}n^2\bar{\phi}^4,$$
and let $\G_{3n}$ denote the corresponding probability measure on $\mathbb{R}^n$. Then the following conclusions under $\G_{3n}$:

\begin{enumerate}
\item[(a)]
 \begin{align}\label{eq:fg2}
 n\E\bar{\phi}^2\lesssim 1,\quad n\E_{\mathbb{G}_{3n}}(\phi_i-\bar{\phi})^2\lesssim 1.
 \end{align}
 
 \item[(b)]
\begin{align}\label{eq:g21}
\sum_{i=1}^n(\phi_i-\bar{\phi})^2\stackrel{P}{\to}4,\\
\label{eq:g22}n^{-1/2}\sum\limits_{1\leq i<j\leq n}(\phi_i+\phi_j-2\bar{\phi})^3\stackrel{P}{\rightarrow}0,\\
\label{eq:g23}\sum\limits_{1\leq i<j\leq n}(\phi_i+\phi_j-2\bar{\phi})^4\stackrel{P}{\to}96.
\end{align}

\item[(c)]
$\sqrt{n}\bar{\phi}\stackrel{D}{\to}\zeta$, where $\zeta$ is a continuous random variable on $\R$ with density proportional to $e^{-\frac{\zeta^2}{2}-\frac{\zeta^4}{24}}$ with respect to Lebesgue measure.

\item[(d)]
$$\sqrt{n}\Big[\sum_{i=1}^n(\phi_i-\bar{\phi})^2-\frac{1}{a_1}\Big]\stackrel{d}{\to}N\Big(0,\frac{2}{a_1^2}\Big).$$

\item[(e)]
For any triangular array $(c_n(1),\ldots,c_n(n))$ with $\sum_{i=1}^nc_{n}(i)=0, \frac{1}{n}\sum_{i=1}^nc_{n}(i)^2\to 1$ we have
$$\sum_{i=1}^nc_{n}(i)\phi_i\stackrel{d}{\to} N\Big(0,\frac{1}{a_1}\Big).$$

\item[(f)]
For every $\varepsilon>0$ we have
$$\limsup_{\delta\to 0}\limsup_{n\to \infty}\P_n(\max_{i,j\in [n]:|i-j|\le n\delta}|S_i(\phi)-S_j(\phi)|>\varepsilon)=0.$$
\end{enumerate}
\end{lem}



Proceeding to verify Lemma \ref{lem:uniqueness_aux} and Lemma \ref{lem:contrast}, we begin by showing the following proposition, which is the analogue of Proposition \ref{ppn:cont1} for $(\theta,\beta)\in \Theta_3$.
\begin{ppn}\label{ppn:cont2}
With $\G_{3n}$ as defined in Lemma \ref{lem:G2} above, we have $\lVert \F_n-\G_{3n}\rVert_{TV}\to 0$.
\end{ppn}
\begin{proof}
Expanding $q(x)=\frac{x^2}{2}-\log\cosh(x)$ by a Taylor's series around $0$ we get \[q(x)=\frac{2x^4}{4!}+R(x),\text { where }|R(x)|\lesssim |x|^6.\] Thus using \eqref{eq:domm} and summing over $1\leq i<j\leq n$ we get
\begin{align}\label{b10}-\log f_n(\phi)=\frac{n}{8}\sum\limits_{i=1}^n(\phi_i-\bar{\phi})^2+\sum_{i<j}\frac{(\phi_i+\phi_j)^4}{2^3 4!}+\sum\limits_{i<j}R(\phi_i+\phi_j).
\end{align} 
Expanding the second term in \eqref{b10} we get
\begin{align*}
\sum_{i<j}(\phi_i+\phi_j)^4=&16N\bar{\phi}^4 +24\sum\limits_{i<j}(\phi_i+\phi_j-2\bar{\phi})^2\bar{\phi}^2\\
+&8\sum_{i<j}(\phi_i+\phi_j-2\bar{\phi})^3\bar{\phi}+\sum_{i<j}(\phi_{i}+\phi_j-2\bar{\phi})^4,
\end{align*}
which along with \eqref{b10} and the identity $\sum_{i<j}(\phi_i+\phi_j-2\bar{\phi})^2=(n-2)\sum_{i=1}^n(\phi_i-\bar{\phi})^2$ gives
\begin{align}\label{b11}
\notag-\log \frac{f_n}{g_{3n}}=&-\frac{n}{24}\bar{\phi}^4+\frac{1}{8}\bar{\phi}^2\Big[(n-2)\sum_{i=1}^n(\phi_i-\bar{\phi})^2-4n\Big]\\
+& \frac{1}{2^34!}\sum_{1\le i<j\le n}\Big[8\bar{\phi}(\phi_i+\phi_j-2\bar{\phi})^3+(\phi_i+\phi_j-2\bar{\phi})^4+R(\phi_i+\phi_j)\Big]
\end{align}
To bound each term on the RHS of \eqref{b11} separately, use Lemma \ref{l31} to get 
\begin{align}
\label{eq:needed_b}&n\bar{\phi}^4\stackrel{P}{\to}0,\quad  \sum_{i<j}R(\phi_i+\phi_j)\lesssim n\sum\limits_{i=1}^n(\phi_i-\bar{\phi})^6+n^2\bar{\phi}^6\to 0,\\
\notag&\Big|(n-2)\bar{\phi}^2\sum_{i=1}^n(\phi_i-\bar{\phi})^2-4n\bar{\phi}^2\Big|\lesssim n\bar{\phi}^2\Big[1+\sum_{i=1}^n(\phi_i-\bar{\phi})^2\Big]=O_P(1),\\
\notag&\Big|\sum_{1\leq i<j\leq n}(\phi_i+\phi_j-2\bar{\phi})^3(2\bar{\phi})\Big|\lesssim \Big|\sqrt{n}\bar{\phi}\Big| \sqrt{n}\sum\limits_{i=1}^n|\phi_i-\bar{\phi}|^3=O_P(1)\\
\notag&\sum\limits_{1\leq i<j\leq n}(\phi_i+\phi_j-2\bar{\phi})^4\lesssim n\sum_{i=1}^n(\phi_i-\bar{\phi})^4=O_P(1).
 \end{align}
 It thus follows that $\log \frac{f_n}{g_{3n}}$ is $O_P(1)$ under $\mathbb{F}_n$. To show the same conclusion under $\G_{3n}$, it suffices to show that the estimates of Lemma \ref{l31} hold under $\G_{3n}$ as well, which follows from part (a) of Lemma \ref{lem:G2}. Thus, using Lemma \ref{ctv} we have that $\F_n$ and $\G_{3n}$ are mutually contiguous.

Finally to show that $\F_n$ and $\G_{3n}$ are close in total variation, invoking (b) of Lemma \ref{ctv} it suffices to show that $\log(f_n/g_{3n})$ converges in probability to a constant under $\G_{3n}$. Invoking \eqref{b11} and \eqref{eq:needed_b}, it suffices to show the following conclusions under $\G_{3n}$:
\begin{align}\label{eq:g21}
\sum_{i=1}^n(\phi_i-\bar{\phi})^2\stackrel{P}{\to}4,\\
\label{eq:g22}n^{-1/2}\sum\limits_{1\leq i<j\leq n}(\phi_i+\phi_j-2\bar{\phi})^3\stackrel{P}{\rightarrow}0,\\
\label{eq:g23}\sum\limits_{1\leq i<j\leq n}(\phi_i+\phi_j-2\bar{\phi})^4\stackrel{P}{\to}96.
\end{align}
But this follows from part (b) of Lemma \ref{lem:G2}. Thus we have verified Proposition \ref{ppn:cont2}.
\end{proof}

 \begin{proof}[Proof of Lemma \ref{lem:uniqueness_aux} and Lemma \ref{lem:contrast} for $(\theta,\beta)\in \Theta_3$]
By Proposition \ref{ppn:cont2} the two probability measures $\F_n$ and $\G_{3n}$ are close in total variation, and so it suffices to work with $\G_{3n}$. But under $\G_{3n}$ the desired conclusions are immediate from parts (c), (d), (e) and (f) of Lemma \ref{lem:G2}. 
 \end{proof}

                                                                                                                                                                                                                                        
\renewcommand\refname{References}            

\bibliographystyle{plain}

\bibliography{bj-sample}            

\section*{Acknowledgements}                 
                                                                     
The first author was partially supported by NSF Grant DMS-1712037.


\begin{appendix}

\section{Proof of Lemma \ref{lem:min} and Proposition \ref{thm:bayes}}\label{sec:four1}

\subsection{Proof of Lemma \ref{lem:min}}
\begin{enumerate}
\item[(a)]
If $(\theta,\beta)\in\Theta_{11}$, then $\beta=0,\theta<1/2$. In this domain, using the inequality $\log\cosh x\leq x^2/2$ we have 
\[q(x)\geq \theta x^2-2\theta^2x^2=\theta x^2(1-2\theta)\geq 0,\] with equality iff $x=0$. Thus $q$ has a unique global minima at $t=0$, and $q''(0)=2\theta(1-2\theta)>0$.
\\If $\theta\in\Theta_{12}$, then $\beta>0$, which gives $q(-t)>q(t)$ for $t>0$. Thus the global minima must lie in $[0,\infty)$. Also since $q(t)$ goes to  $\infty$ as $t\rightarrow\infty$, the global minima is not attained at $\infty$. If there is a local minima at $t$ for some $t> 0$, then it must satisfy  
$t=\tanh[2\theta t+\beta]$ which has a unique strictly positive solution on $(0,\infty)$. Finally note  that $q'''(t)\neq 0$, and so for $t$ to be a minima we must have that $q''(t)>0$.
This completes the proof for $\theta\in\Theta_{12}$.

\item[(b)]
For $(\theta,\beta)\in \Theta_2$, differentiating $q$ we get $q'(x)=2\theta[x-2\tanh(\theta x)]$ which has  exactly three real roots $0,\pm t$ where $t$ is a root of $t=\tanh(2\theta t)$. Also note that $q''(0)<0$, and so $\pm t$ are local minima of $q(.)$ and $0$ is a local maxima. Since $q(t)\rightarrow \infty$ as $|t|\rightarrow \infty$, we have that $\pm t$ are also the global minima, as claimed. A similar argument as above then shows that $q''(t)>0$.

\item[(c)]
If  $(\theta,\beta)\in \Theta_{3}$, then $\beta=0,\theta=1/2$, and so $q(x)\ge \theta x^2-2\theta^2 x^2=0$, with equality iff $x=0$. Thus $q(.)$ has a unique global minimum at $0$, with $q'(0)=q''(0)=0$.
\end{enumerate}

\subsection{Proof of Proposition \ref{thm:bayes}}

\begin{enumerate}
\item[(a)]
Using \eqref{Y:edge_star} and \eqref{eq:phi_z}, the joint likelihood of $(Y,\phi)$ is proportional to
\begin{align}
\notag&\exp\left\{\frac{\theta}{2(n-1)}\sum_{i=1}^nk_i(y)^2+\frac{\beta}{2}\sum_{i=1}^nk_i(y)-\frac{(n-1)\theta}{2}\sum_{i=1}^n\Big(\phi_i-\frac{k_i(y)}{n-1}\Big)^2\right\}\\
\notag=&\exp\left\{ -\frac{(n-1)\theta}{2}\sum_{i=1}^n\phi_i^2+\sum_{i=1}^nk_i(y)\Big(\frac{\beta}{2}+\phi_i\Big)\right\}\\
\label{eq:joint2020}=&\exp\left\{-\frac{(n-1)\theta}{2}\sum_{i=1}^n\phi_i^2+\sum_{i<j}y_{ij}(\phi_i+\phi_j+\beta)\right\}.
\end{align}
It follows from \eqref{eq:joint2020} that given $\phi$ the random variables $\{Y_{ij}\}_{1\le i<j\le n}$ are mutually independent, and have the conditional distribution as specified in part (a).

\item[(b)]
Using \eqref{eq:joint2020}, the marginal density of $\phi$ is proportional to
\begin{align*}
&\sum_{y\in \{-1,1\}^{n\choose 2}}\exp\left\{-\frac{(n-1)\theta}{2}\sum_{i=1}^n\phi_i^2+\sum_{i<j}y_{ij}(\phi_i+\phi_j+\beta)\right\}\\
=&\exp\left\{-\frac{(n-1)\theta}{2}\sum_{i=1}^n\phi_i^2\right\}\prod_{1\le i<j\le n}\sum_{y_{ij}\in \{-1,1\}} \exp\left\{ y_{ij}(\phi_i+\phi_j+\beta)\right\}\\
=&\exp\left\{-\frac{(n-1)\theta}{2}\sum_{i=1}^n\phi_i^2\right\}2^{n\choose 2}  \exp\left\{\sum_{i<j} \log\cosh\Big(\theta(\phi_i+\phi_j)+\beta\Big)\right\}=2^{n\choose 2} f_n(\phi),
\end{align*}
which verifies part (b).
\end{enumerate}

\section{Proof of Lemma \ref{ctv}, Lemma \ref{lem:G1} and Lemma \ref{lem:12whole}}\label{sec:four2}

\subsection{ {Proof of Lemma \ref{ctv}}}

\begin{enumerate}
\item[(a)]
Let $A_n$ be a sequence of sets such that $\mathbb{G}_n(A_n)\rightarrow 0$. We will show that $\mathbb{H}_n(A_n)\rightarrow 0$, which will give $\mathbb{H}_n\ll \mathbb{G}_n$. The other implication then follows by symmetry. Fix $\epsilon\in (0,1)$ arbitrary, and let $M=M(\epsilon)$ be such that setting
$$ B_n(\epsilon):=\left\{e^{-M}\leq L_n\leq e^M\geq 1-\epsilon\right\}\text{ with }L_n:=\frac{h_n}{g_n}.$$
Using the fact that $\log L_n$ is $O_P(1)$ under both probability measures, 
we have $\mathbb{H}_n(B_n)\geq 1-\epsilon$, and $\mathbb{G}_n(B_n)\geq 1-\epsilon.$ 
Thus setting $\epsilon=1/2$ and $b_n:=\int g_n d\mu$ and $a_n=\int h_n d\mu$ we have
\begin{align*}
b_n=\int g_nd\mu_n\leq e^{M(\frac{1}{2})}\int_{B_n(\frac{1}{2})}h_nd\mu_n+ \frac{b_n}{2}\leq e^{M(\frac{1}{2})} a_n+ \frac{b_n}{2}.
\end{align*}
This gives $b_n\leq 2e^{M(1/2)}a_n$, i.e. $b_n/a_n$ is bounded above. Now for any $\epsilon>0$ we have

\begin{align*}
\mathbb{H}_n(A_n)\leq &\frac{\int_{A_n\cap B_n(\epsilon)}h_nd\mu_n}{a_n}+\epsilon\leq e^{M(\epsilon)} \mathbb{G}_n(A_n)\frac{b_n}{a_n}+\epsilon\leq2e^{M(\epsilon)+M(1/2)}\mathbb{G}_n(A_n)+\epsilon
\end{align*}
Taking $\limsup$ on both sides as $n\rightarrow\infty$, we have $\limsup_{n\to\infty} \mathbb{H}_n(A_n)\leq \epsilon$. Since $\epsilon<1$ is arbitrary, we have the result.

\item[(b)]
Note that the above proof implies that $r_n:=\log a_n-\log b_n$ is a bounded sequence of reals, and so converges along a  subsequence to $r$, say. Without loss of generality we restrict to this subsequence, and set \[\overline{h}_n:=\frac{h_n}{a_n},\quad \overline{g}_n:=\frac{g_n}{b_n},\quad  \overline{L}_n:=\log \overline{h}_n-\log \overline{g}_n=L_n-r_n.\]  Then we have 
\[(X_{n},\overline{L}_{n})\stackrel{d,\mathbb{G}_{n}}{\rightarrow}N(\mu_1,\mu_2-r,\sigma_1^2,\sigma_2^2,\sigma_{12})\] But then \cite{roussas1972contiguity} gives $r=\mu_2+\frac{\sigma_2^2}{2}$, and so $r_n$ converges to $\mu_2+\frac{\sigma_2^2}{2}$. But then \[(X_{n},\overline{L}_{n})\stackrel{d,\mathbb{G}_{n}}{\rightarrow}N(\mu_1,\mu_2-r,\sigma_1^2,\sigma_2^2,\sigma_{12})\] from which we have by \cite{roussas1972contiguity} we get \[X_n\stackrel{d,\mathbb{H}_n}{\rightarrow}N(\mu_1+\sigma_{12},\sigma_1^2).\]

\item[(c)]
Since $\mathbb{H}_n$ and $\mathbb{G}_n$ are mutually contiguous, it follows that $L_n\stackrel{d,\mathbb{H}_n}{\rightarrow}c$. Fix $\epsilon,\delta>0$, arbitrary. Then for all large $n$ we have 
$$\mathbb{H}_n(B_n)>1-\epsilon,\mathbb{G}_n(B_n)>1-\epsilon,\quad B_n:=\Big\{c-\delta<\log\Big(\frac{h_n}{g_n}\Big)<c+\delta
\Big\}.$$ 
Thus
\begin{align*}
b_n=&\int_{B_n}g_nd\mu_n+\int_{B_n^c}g_nd\mu_n
\leq e^{-c+\delta}\int_{B_n}h_nd\mu_n+\epsilon b_n
\leq e^{-c+\delta}a_n+\epsilon b_n
\end{align*}
which gives  $\frac{b_n}{a_n}\leq \frac{e^{-c+\delta}}{1-\epsilon}$.

By similar calculations, we have $\frac{a_n}{b_n}\leq \frac{e^{c+\delta}}{1-\epsilon}$. Since $\delta,\epsilon>0$ arbitrary, we have $\frac{a_n}{b_n}\rightarrow e^c$, which gives ${\overline{h}_n}/{\overline{g}_n}\stackrel{d}{\rightarrow }1$ by Slutsky, under both $\mathbb{H}_n$ and $\mathbb{G}_n$. The desired conclusion is immediate from this.
\end{enumerate}


\subsection{Proof of Lemma \ref{lem:G1} }

Since $-\log g_{1n}(.)$ is quadratic, 
we have $\phi-t:=(\phi_1-t,\ldots,\phi_n-t)\sim N({\bf 0},\Gamma),$ with $\Gamma^{-1}:=na_1{\bf I}-a_2 {\bf 1} {\bf 1}'$. Inverting, we have $\Gamma=\frac{1}{na_1}{\bf I}+\frac{a_2}{n^2a_1(a_1-a_2)}$, and so 
\begin{align}\label{eq:rep}
(\phi_1-t,\ldots,\phi_n-t)\stackrel{d}{=}\frac{1}{\sqrt{na_1}}(W_1,\ldots, W_n)+\sqrt{\frac{a_2}{n^2a_1(a_1-a_2)}}W_0,
\end{align} 
where $(W_i)_{0\le i\le n}\stackrel{i.i.d.}{\sim}N(0,1)$. Note that we have tacitly used the fact that $a_1>a_2$, which follows from Lemma \ref{lem:min} along with the observation that $$q''(t)=2\theta-4\theta^2(1-t^2)=2(a_1-a_2).$$

\begin{enumerate}
\item[(a)]
Using \eqref{eq:rep} we can write
\begin{align}
\label{eq:rep3}n(\bar{\phi}-t)=&\sqrt{\frac{1}{na_{1}}}\sum\limits_{i=1}^{n}W_{i}+\sqrt{\frac{a_{2}}{a_1(a_1-a_2)}}W_0,\\
\label{eq:rep2}\sqrt{n}[\sum\limits_{i=1}^{n}(\phi_{i}-\bar{\phi})^{2}-a_{1}^{-1}]=&\frac{1}{a_{1}\sqrt{n}}\sum\limits_{i=1}^{n}(W_{i}^{2}-1)-\frac{(\sqrt{n}\bar{W})^{2}}{\sqrt{n}a_{1}},
\end{align}
where $\bar{W}:=\frac{1}{n}\sum_{i=1}^nW_i$.
The desired conclusions of part (a) are immediate from these representations.

\item[(b)]
Again using \eqref{eq:rep} gives
\begin{align*}
 \sum_{1\le i<j\le n}\Big(\phi_i+\phi_j-2t\Big)^4 &=\sum\limits_{i<j}\left[\sqrt{\frac{1}{na_{1}}}W_{i}+\sqrt{\frac{1}{na_{1}}}W_{j}+2\sqrt{\frac{a_{2}}{n^{2}a_{1}(a_1-a_{2})}}W_0\right]^{4}    \\
 &=\frac{1}{n^2a_1^2}\sum_{1\le i<j\le n}(W_i+W_j)^4+o_P(1)\stackrel{P}{\to}\frac{6}{a_1^2}.
\end{align*}

\item[(c)]
Again using the representation \eqref{eq:rep} gives
\begin{align}\label{eq:rep1}
\notag n\sum\limits_{i=1}^{n}(\phi_{i}-t)^{3}=&\frac{1}{\sqrt{na_{1}^{3}}}\sum\limits_{i=1}^{n}W_{i}^{3}+\frac{3\sqrt{a_{2}}W_0\sum\limits_{i=1}^{n}W_{i}^{2}}{na_{1}\sqrt{a_1(a_1-a_{2})}}+O(n^{-\frac{3}{2}})W_0^{2}\sum\limits_{i=1}^{n}W_{i}+O(n^{-2})W_0^{3}\\
=&\frac{1}{\sqrt{na_{1}^{3}}}\sum\limits_{i=1}^{n}W_{i}^{3}+\frac{3}{a_{1}}\sqrt{\frac{a_{2}}{a_{1}(a_1-a_{2})}}W_0+o_{p}(1),
\end{align}
and 
\begin{align}\label{eq:rep5}
\sum_{i=1}^nc_n(i)\phi_i=\frac{1}{\sqrt{na_1}}\sum_{i=1}^nc_n(i)W_i.
\end{align}
Combining \eqref{eq:rep3}, \eqref{eq:rep2}, \eqref{eq:rep1} and \eqref{eq:rep5} gives
\small\begin{align*}
&\left[n(\bar{\phi}-m),\sqrt{n}\Big(\sum\limits_{i=1}^{n}(\phi_{i}-\bar{\phi})^{2}-a_{1}^{-1}\Big),n\sum\limits_{i=1}^{n}(\phi_{i}-m)^{3},\sum_{i=1}^nc_n(i)\phi_i\right]\\
=&\bigg[ \sqrt{\frac{1}{na_{1}}}\sum\limits_{i=1}^{n}W_{i}+\sqrt{\frac{a_{2}}{a_{1}(a_1-a_{2})}}W_0,\frac{1}{a_{1}\sqrt{n}}\sum\limits_{i=1}^{n}(W_{i}^{2}-1), \\
& \frac{1}{\sqrt{na_{1}^{3}}}\sum\limits_{i=1}^{n}W_{i}^{3}+\frac{3}{a_{1}}\sqrt{\frac{a_{2}}{a_{1}(a_1-a_{2})}}W_0 ,\frac{1}{\sqrt{na_1}}\sum_{i=1}^nc_n(i)W_i\bigg]+o_p(1).
\end{align*}
The desired conclusion follows from this on invoking the Central Limit Theorem to note that
\begin{equation*}
\left[\sqrt{\frac{1}{na_{1}}}\sum\limits_{i=1}^{n}W_{i},\sqrt{\frac{1}{na_{1}^{2}}}\sum\limits_{i=1}^{n}(W_{i}^{2}-1),  \sqrt{\frac{1}{na_{1}^{3}}}\sum\limits_{i=1}^{n}W_{i}^{3},\frac{1}{\sqrt{n}}\sum_{i=1}^nc_n(i)W_i  \right]\stackrel{D}{\rightarrow}N(0,\widetilde{\Sigma})
\end{equation*}
with
%
\begin{equation*}\widetilde{\Sigma}:=
 \left[ {\begin{array}{cccc}
   \frac{1}{a_1}& 0 & \frac{3}{a_1^2}&0 \\
   0& \frac{2}{a_1^2} & 0&0 \\
   \frac{3}{a_1^2}&0&\frac{15}{a_1^3}&0\\
   0&0&0&\frac{1}{a_1}
  \end{array} } \right],
\end{equation*}
along with continuous mapping theorem.

\item[(d)]

Using \eqref{eq:rep} we have
$$\max_{i,j\in [n]:|i-j|\le n\delta}|S_i(\phi)-S_j(\phi)|=\frac{1}{\sqrt{na_1}} \max_{i,j\in [n]}|\tilde{S}_i({\bf W})-\tilde{S}_j({\bf W})|,$$
where ${\bf W}:=(W_1,\ldots,W_n)$, and $\tilde{S}_i({\bf W}):=\sum_{j=1}^iW_j$. This gives
$$\P_n(\max_{i,j\in [n]:|i-j|\le n\delta}|S_i(\phi)-S_j(\phi)|>\varepsilon)\le \P_n( \max_{i,j\in [n]}|\tilde{S}_i({\bf W})-\tilde{S}_j({\bf W})|>\varepsilon \sqrt{na_1}),$$
where the RHS above goes to $0$ under the double limit $n\to\infty$ followed by $\delta\to 0$, invoking tightness of sample paths of i.i.d.~sums as in the proof of Theorem \ref{bm}.
\end{enumerate}

\subsection{Proof of Lemma \ref{lem:12whole} }

For proving Lemma \ref{lem:12whole} we need a second moment bound for the conditionally centered sum of $\{Y_e\}_{e\in \mathcal{E}}$, where $\mathcal{E}$ is the collection of all $(i,j)$ with $1\le i<j\le n$. 
\begin{defn}
For $e=(i,j)\in \mathcal{E}$,  let $N(e):=\{(i,k), (j,k), k\ne i,j\}$ denote the collection of pairs $(a,b)\in \mathcal{E}$ with exactly one element from $\{i,j\}$, and let \[t_{e}(y):=\frac{1}{2(n-1)}\sum_{k\ne i,j}\{y_{ik}+y_{jk}\}=\frac{k_i({ y})+k_j({y})-2y_{ij}}{2(n-1)}.\]
\end{defn}
\begin{lem}\label{pair}
For any $(\theta,\beta)$ we have
  \[\E\Big[\sum\limits_{e\in \mathcal{E}}\Big\{Y_{e}-\tanh(2\theta t_e+\beta)\Big\}\Big]^{2}=O(n^2).\]
   \end{lem}
The proof of Lemma \ref{lem:12whole} is deferred to the end of this section.

\begin{proof}[Proof of Lemma \ref{lem:12whole}]
\begin{enumerate}
\item[(a)]
To begin, note that \[p(x,y)=\frac{\theta}{4}(x-y)^2+q\Big(\frac{x+y}{2}\Big),\text{ with }q(x)=\theta x^2-\log\cosh(2\theta x+\beta)\] as defined in Lemma \ref{lem:min}.
Now the function $\tilde{q}(.):\R\mapsto \R$ defined by 
\[\tilde{q}(x):=\frac{q(x)-q(t)}{(x-t)^2},x\ne t,\quad \tilde{q}(t):=\frac{1}{2}q''(t)\]
is a continuous function which satisfies $\lim_{x\rightarrow\pm \infty}\tilde{q}(x)=\frac{\theta}{2}$. Also $\tilde{q}(.)$ is strictly positive everywhere on $U$ by Lemma \ref{lem:min}. Thus there exists constants $\tilde{\lambda}_1\ge \tilde{\lambda}_2$ such that for all $x\in U$ we have
$\tilde{\lambda}_2(x-t)^2\le q(x)\le \tilde{\lambda}_1(x-t)^2$, which in turn shows that for all $(x,y)\in U^2$ we have
\begin{align*}
p(x,y)\le \frac{\theta}{4}(x-y)^2+\frac{\tilde{\lambda}_1}{4}(x+y-2t)^2\le \frac{\max(\theta,\tilde{\lambda}_1)}{2}\Big[(x-t)^2+(y-t)^2\Big],\\
p(x,y)\ge \frac{\theta}{4}(x-y)^2+\frac{\tilde{\lambda}_2}{4}(x+y-2t)^2\ge \frac{\min(\theta,\tilde{\lambda}_2)}{2}\Big[(x-t)^2+(y-t)^2\Big].
\end{align*}
This verifies part (a).
\\

\item[(b)]

Use part (a) to note that for all $\phi\in U^n$ we have 
   \[\frac{(n-1)\lambda_2}{2}\sum\limits_{i=1}^n(\phi_i-t)^2\leq -\log f_n(\phi)  \leq \frac{(n-1)\lambda_1}{2}\sum\limits_{i=1}^n(\phi_i-t)^2,\] and so for  any $M>0$,
\begin{align}
\notag\mathbb{P}_{n,U}\Big(\sum\limits_{i=1}^n(\phi_i-t)^2>M\Big)= &\frac{\int_{U^n}e^{-f_n(\phi)}\mathbf{1}\{\sum_{i=1}^n(\phi_i-t)^2\}d\phi}{\int_{U^n}e^{-f_n(\phi)}d\phi}\\
\notag\leq&\frac{\int_{U^n}e^{-\frac{(n-1)\lambda_2}{2}\sum\limits_{i=1}^n(\phi_i-t)^2}\mathbf{1}\{\sum_{i=1}^n(\phi_i-t)^2\}d\phi}{\int_{U^n}e^{-\frac{(n-1)\lambda_1}{2}\sum\limits_{i=1}^n(\phi_i-t)^2}d\phi}\\
\label{eq:step0}\leq &\Big(\frac{\lambda_1}{\lambda_2}\Big)^{\frac{n}{2}}\frac{\mathbb{P}_n\Big((\chi_n^2\geq (n-1)\lambda_2 M\Big)}{\P_n\Big(t+\frac{W_0}{\sqrt{(n-1)\lambda_1}}\in U\Big)^n}
\end{align} 
where $W_0\sim N(0,1)$ and $\chi_n^2$ is a chi-square random variable with $n$ degrees of freedom. Since $t\in U$, using standard tail estimates of $N(0,1)$ for $(\theta,\beta)\in \Theta_2$ we have
\begin{align}\label{eq:denom1}
\P_n\left(t+\frac{W_0}{\sqrt{(n-1)\lambda_1}}\in U\right)^n\ge \P_n(|W_0|\le t\sqrt{(n-1)\lambda_1})^n\rightarrow 1.
\end{align}
and so the denominator in the RHS of \eqref{eq:step0} converges to $1$.  
%
Proceeding to bound the numerator of \eqref{eq:step0}, use the moment generating function of $\chi_n^2$ along with  Markov's inequality  to get
\[\limsup\limits_{n\rightarrow\infty}\frac{1}{n}\log \mathbb{P}_{n,U}(\sum\limits_{i=1}^n(\phi_i-\phi_0)^2>M)\leq \log\Big(\frac{\lambda_1}{\lambda_2}\Big)-\frac{\lambda_2 M}{2}+\frac{1}{2}\log(\lambda_2 M).\] Since the RHS above converges to $-\infty$ as $M\rightarrow\infty$, there exists $M$ depending on $(\lambda_1,\lambda_2)$ such that the RHS above is  negative,  from which the conclusion of part (b) follows. 
\\

\item[(c)]
A direct calculation gives
\begin{align}\label{eq:step1}
\notag&\mathbb{E}_{n,U}|\phi_1-t|^l\\
\notag=&\mathbb{E}_{n,U}|\phi_1-t|^l\mathbf{1}\{\sum\limits_{i=2}^n(\phi_i-\phi_t)^2>M\}+\mathbb{E}_{n,U}|\phi_1-t|^l\mathbf{1}\Big\{\sum\limits_{i=2}^n(\phi_i-t)^2\leq M\Big\}\\
\le &\sqrt{\mathbb{E}_{n,U}(\phi_1-t)^{2l}}\sqrt{\mathbb{P}_{n,U}(\sum\limits_{i=2}^n(\phi_i-t)^2>M)}+\mathbb{E}_{n,U}|\phi_1-\phi_0|^l\mathbf{1}\Big\{\sum\limits_{i=2}^n(\phi_i-t)^2\leq M\Big\} ,   
\end{align}
where the last step uses Cauchy-Schwarz inequality. To bound the first term in the RHS of \eqref{eq:step1},  with $W_0\sim N(0,1)$ we have \[|\phi_1-t|\stackrel{d}{=}\Big|\frac{k_i}{n-1}-t+\frac{W_0}{\sqrt{(n-1)\theta}}\Big|\le 1+\frac{|W_0|}{\sqrt{(n-1)\theta}}.\] Since $\P_n(\phi\in U)\ge \frac{1}{2}$, it follows from the above display that 
\begin{align}\label{eq:crude}
\P_{n,U}(|\phi_1-t|>2)\le e^{-\Omega(n)}, \quad \E_{n,U}(\phi_1-t)^{2l}\lesssim_l 1.
\end{align}
 Using part (b) then implies that the first term in the RHS of \eqref{eq:step1} is bounded by $e^{-\Omega(n)}$. For estimating the second term, note that
 the conditional density of $(\phi_1|\phi_i,i\geq 2)$ is proportional to $\prod\limits_{i=2}^ne^{-p(\phi_1,\phi_i)}$. Applying part (a) gives
\begin{align*}
\mathbb{E}_{n,U}\Big(|\phi_1-t|^l\Big|\phi_i,i\geq 2\Big)
 =&\frac{\int_{U}|\phi_1-t|^l\prod\limits_{i=2}^ne^{-p(\phi_1,\phi_i)}d\phi_1}{\int_{U}\prod\limits_{i=2}^ne^{-p(\phi_1,\phi_i))}d\phi_1}\\
  \leq &e^{\frac{\lambda_1-\lambda_2}{2}\sum\limits_{i=2}^n(\phi_i-t)^2}\sqrt{\Big(\frac{\lambda_1}{\lambda_2}\Big)}\frac{\mathbb{E}|W_0|^l}{[(n-1)\theta]^{l/2}\mathbb{P}\left(\frac{|W_0|}{\sqrt{(n-1)\lambda_1}}\in U\right)}
 \end{align*}
 where  $W_0\sim N(0,1)$. Since we have $\sum\limits_{i=2}^n(\phi_i-\phi_0)^2\leq M$ on the conditioned set and $\mathbb{P}_n(\frac{|W_0|}{\sqrt{(n-1)\lambda_1}}\in U)$ converges to $1$ using \eqref{eq:denom1},  the conclusion of part (c) follows.
\\

 
%

\item[(d)]
Expanding $\tanh(.)$ by Taylor's series gives
\begin{align*}
\tanh(2\theta t_e(y)+\beta)=\tanh(2\theta t+\beta)+2\theta (t_e-t)\text{sech}^2(2\theta t+\beta)+O(|t_e-t|^2)\\
=t+2\theta(1-t^2)(t_e-t)+O(|t_e-t|^2).
\end{align*}
This gives
\begin{align}\label{u1}
\Big|\sum_{e\in E}(y_e-\tanh(2\theta t_e(y)+\beta))-
\sum\limits_{e\in E}(y_e-t)-\sum\limits_{e\in E}2\theta(1-t^2)\Big(\frac{t_e(y)}{2(n-1)}-t\Big)\Big|\lesssim \sum_{e\in E}(t_e-t)^2.
\end{align} 
 If $e=(i,j)$, then we have $$t_e(y)=\frac{k_i+k_j-2y_{ij}}{2(n-1)}=\frac{\phi_i+\phi_j}{2}-\frac{W_i+W_j}{2\sqrt{(n-1)\theta}}-\frac{y_{ij}}{n-1},$$ and so the RHS of \eqref{u1} is bounded (up to constants) by
\begin{align}\label{u2}
n\Big[\sum\limits_{i=1}^n(\phi_i-t)^2+\sum\limits_{i=1}^n\frac{W_i^2}{(n-1)\theta}+1\Big]=:nV_n,\quad \E_{n,U}V_n^2\lesssim 1,
\end{align}
where the second conclusion above uses part (a). 
Since the LHS of \eqref{u1} equals 
\begin{align*}\Big|\sum\limits_{e\in E}(y_e-\tanh(2\theta t_e+\beta)-(1-2\theta(1-t^2))\sum\limits_{e\in E}(y_e-t)-\frac{2\theta(1-t^2)}{n-1}\sum\limits_{e\in E}y_e\Big)\Big|,
\end{align*}
it follows on invoking \eqref{u1} and \eqref{u2} along with Lemma \ref{pair} that 
\begin{align*}
\E_{n,U}\big(\sum\limits_{e\in E}(y_e-t)\big)^{2}\lesssim n^{2}. 
\end{align*} from which the desired conclusion is immediate. 
\end{enumerate}
\end{proof}

\subsection{Proof of Lemma \ref{pair}}

Letting \[H_n(y)=\frac{1}{N}\sum_{e\in E}(y_e-\tanh\Big(\frac{\theta}{n-1}t_e(y)+\beta\Big)\]  it suffices to show that $\mathbb{E}H_n(Y)^2\lesssim \frac{1}{N}$. The technique for the proof is by using exchangeable pairs, and is adapted from \cite{chatterjee2007estimation}.

Produce an exchangeable pair $(Y,Y')$ in the following  way: 
\\
 Let $Y\sim \mathbb{P}_{n}$. To simulate $Y'$, choose an $e\in {\mathcal E}$ uniformly at random, and replace $Y_e$ by $Y_e'$, which is a simulation from the conditional distribution of $(Y_e|Y_f,f\neq e,f\in {\mathcal E})$.  Let $F(Y,Y')=E(Y)-E(Y')$. Note that \[\mathbb{P}_n(Y_e=1|Y_f=y_f,f\neq e)=\frac{e^{2t_e(y)\theta+\beta}}{e^{2\theta t_e(y)+\beta}+e^{-2\theta t_e(y)-\beta}}\] and so
\begin{align*}
\mathbb{E}[F(Y,Y')|Y=y]=\mathbb{E}(E(Y)-E(Y')|Y=y)=\frac{1}{N}\sum_{e\in E}\big{[}y_e-\tanh(2\theta t_e(y)+\beta)\big{]}=H_n(y).
\end{align*}
Also note that
\begin{align*}
\E H_n(Y)^2=\mathbb{E}H_n(Y)F(Y,Y')=\mathbb{E}H_n(Y')F(Y',Y)=-\mathbb{E}H_n(Y')F(Y,Y'),
\end{align*}
where we use the fact that $(Y,Y')$ are exchangeable and $F$ is antisymmetric. This readily implies  \begin{align}\label{fg1}\E H_n(Y)^2=\frac{1}{2}\mathbb{E}[(H_n(Y)-H_n(Y'))F(Y,Y')]=:\mathbb{E}v_n(Y)\end{align} with 
$v_n(y):=\frac{1}{2}\mathbb{E}[(H_n(Y)-H_n(Y'))F(Y,Y')|Y=y]$. Proceeding to estimate 
 this,  let $y^e$ denote  $y$ with the sign of $y_e$ reversed, and let \[p_e(y):=\frac{ e^{-y_e[2\theta t_e(y)+\beta]}}{e^{2\theta t_e(y)+\beta}+e^{-2\theta t_e(y)-\beta}}=\P_n(Y_e=-y_e|Y_f=y_f,f\neq e).\] Then we have
\begin{align*}
2v_n(y)=\frac{1}{N}\sum\limits_{e\in {\mathcal E}}(H_n(y)-H_n(y^e))F(y,y^e)p_e(y)=\frac{1}{N}\sum\limits_{e\in {\mathcal E}}(H_n(y)-H_n(y^e))2y_ep_e(y),
\end{align*}
where
\begin{align*}
H_n(y)-H_n(y^e)=\frac{2y_e}{N}+\frac{1}{N}\sum\limits_{f\in N(e)}\left[-\tanh(2\theta t_f(y)+\beta)+\tanh(2\theta t_f(y^e)+\beta)\right]
\end{align*}
The first term in the display above is bounded by $\frac{2}{N}$. Also, using $|\tanh(a)-\tanh(b)|\leq |a-b|$ for all $a,b\in \mathbb{R}$, it follows that the second term is bounded by $\frac{2|\theta|}{N}$. Thus we have $|v_n(y)|\leq \frac{1+|\theta|}{N}$, which along with (\ref{fg1})  completes the proof of the lemma.





\section{Proof of Proposition \ref{lem:sep} }\label{sec:four3}

\subsection{Proof of Proposition \ref{lem:sep}}
Using the identity \eqref{eq:domm} along with part (b) of Lemma \ref{lem:min}, it follows that $p(x,y)$ has two global minima at $\pm(t,t)$. We now break the proof into steps.
\\

\begin{itemize}
\item[{\bf Step 1}]
Letting $I=[t/2,3t/2]$, we first show that there exists a positive integer $M<\infty$ such that 
\begin{align}\label{fg3}
\mathbb{P}_{n}(A_n)\leq e^{-\Omega(n)},\quad A_n:=\Big\{\sum_{i=1}^n 1\{|\phi_i|\notin I\}>M\Big\}.
\end{align}
This means that with very high probability atmost finitely many of the co-ordinates of $\phi$ are in $\pm I$. 

{\bf Proof:}
For proving \eqref{fg3}, 
 note that $p(x,y)\geq p(|x|,|y|)$ for any $x,y\in \R^2$, and so for any quadrant $\mathcal{Q}$ we have
\begin{align*}
\int\limits_{\mathcal{Q}\cap A_n}f_n(\phi)d\phi\leq \int\limits_{\mathcal{Q}_1\cap {A}_n}f_n(\phi)d\phi,
\end{align*}
where $\mathcal{Q}_1$ is the first quadrant in $\mathbb{R}^n$. 
Concentrating on $x,y>0$, since the unique global minima of $p(x,y)$  is at $(t,t)$, using arguments similar to part (a) of Lemma \ref{lem:12whole} gives the existence of finite positive constants $\lambda_1, \lambda_2$ such that for all $x,y\ge 0$ we have 
\[\frac{\lambda_2}{2}[(x-t)^2+(y-t)^2\leq p(x,y)\leq \frac{\lambda_1}{2}[(x-t)^2+(y-t)^2].\] This gives
\begin{align*}
\frac{\int\limits_{\mathcal{Q}_1\cap A_n}f_n(\phi)d\phi}{\int\limits_{I^n}f_n(\phi)d\phi}\leq \Big(\frac{\lambda_1}{\lambda_2}\Big)^{n/2}\frac{\mathbb{P}\left(\sum_{i=1}^n1\Big\{\frac{W_i}{\sqrt{(n-1)\lambda_2}}\Big\}>t/2\}>M \right)}{\mathbb{P}_n\Big((\frac{|W_0|}{\sqrt{(n-1)\lambda_1}}<t/2\Big)^n}
\end{align*}
where $\{W_i\}_{0\le i\le n}\stackrel{i.i.d.}{\sim}N(0,1)$. The probability in the denominator converges to $1$ as before. By a union bound, the probability in the numerator is bounded by $$\Big({{n}\atop{M}}\Big)e^{-M\Omega(n)}.$$ Since there are $2^n-2$ co-ordinates other than $\pm\mathcal{Q}_1$, we have 
\[\mathbb{P}_{n}\Big(\sum_{i=1}^n\{|\phi_i|\notin I,1\leq i\leq n\}>M)\leq 2^n\Big(\frac{\lambda_1}{\lambda_2}\Big)^{n/2}\Big({{n}\atop{M}}\Big)e^{-M\Omega(n)}.\]
Choosing $M$ large enough gives \eqref{fg3}, as desired.
\\

\item[{\bf Step 2}]
Letting $I_1=I_1(\phi):=\{i:\phi_i\in I\}$, $I_2=I_2(\phi):=\{i:\phi_i\in -I\}$, we now show that there exists $N<\infty$ such that
\begin{align}\label{fg4}
\mathbb{P}_{n}(I_1>N,I_2>N)\leq e^{-\Omega(n)}.
\end{align} 

{\bf Proof:}
To show \eqref{fg4} first note that   there exists $c>0$ such that  
\begin{align}\label{dd1}
\min\Big(\inf_{x\in I,y\in -I}\{p(x,y)-p(|x|,|y|)\},\inf_{x\in -I,y\in I}\{p(x,y)-p(|x|,|y|)\}\Big)\geq c.
\end{align}
Indeed, (\ref{dd1}) follows from the fact that $p(x,y)>p(|x|,|y|)$ for all $(x,y)$ on $\pm\{\overline{I}\times -\overline{I}\}$  which is compact, and  $p(x,y)> p(|x|,|y|)$ for all $(x,y)$ in this domain. Using \eqref{dd1}  we have 
$f_n(\phi)\leq f_n(|\phi|)e^{-I_1I_2c}$. Now
\begin{align}\label{fg5}
\mathbb{P}_{n}(I_1>N,I_2>N)\leq&\mathbb{P}_{n}(I_1>N,I_2>N,A_n^c)+\mathbb{P}_{n}(A_n),
\end{align}
with the second term  bounded by $e^{-\Omega(n)}$ by \eqref{fg3}. For the first term note that if $I_1>N,I_2>N$ and $I_1+I_2\geq n-M$, then $I_1I_2>N(n-M-N)$, and so 
\begin{align*}
\mathbb{P}_{n}(I_1>N,I_2>N,A_n^c)\leq \frac{2^n  e^{-N(n-M-N)c}\int\limits_{\mathcal{Q}_1}f_n(\phi)d\phi}{\int\limits_{\mathcal{Q}_1}f_n(\phi)d\phi}=2^ne^{-N(n-M-N)c}.
\end{align*}
Thus choosing $N$ large enough gives \eqref{fg4}). 
\\

\item[{\bf Step 3}]
Setting $J_1:=\{i:\phi_i>0\},J_2=\{i:\phi_i<0\}$, we will show that
\[\P_n(J_1<n,J_2<n)\leq e^{-\Omega(n)},\] which will prove the Lemma.

{\bf Proof:}
Combining \eqref{fg3} and \eqref{fg4} readily gives
\begin{align}\label{fg6}
\mathbb{P}_{n}(I_1<n-R,I_2<n-R)\leq e^{-\Omega(n)},
\end{align}
with $R:=M+N$, i.e. with high probability at least $n-R$ of the $\phi_i$'s are in one of $\pm I$. Since
$J_1\geq I_1,J_2\geq I_2$, \eqref{fg6} gives
\begin{align*}
\mathbb{P}_{n}(J_1<n,J_2<n)\leq2\mathbb{P}_{n}(J_1<n,J_2<n,I_1\geq n-R)+\mathbb{P}_{n}(I_1<n-R,I_2<n-R)
\end{align*}
The second term is $e^{-\Omega(n)}$ by \eqref{fg5}. Turning to deal with the first term, note that for $x\in I,y\leq 0$ we have 
\begin{align}\label{dd2}\inf_{x\in I,y\leq 0}p(|x|,|y|)-p(x,y)\geq \tilde{c}>0.\end{align}
 Indeed, as before this function is positive point-wise on compact subsets of $I\times (-\infty,0)$, and their difference goes to $\infty$ if $y\rightarrow-\infty$. 
\\

Thus, since $J_1+J_2=n$, and $I_1>n-R$, we have that $J_2\geq  1$, and so on this set there exists at least $(n-1)$ pairs $(i,j)$ such that $\phi_i\in I,\phi_j<0$. This readily implies by \eqref{dd2} that $f_n(\phi)\leq f_n(|\phi|)e^{-(n-1)\tilde{d}}$. Also note that $I_1\geq n-R$ can occur only on at most $2^{R}\Big({{n}\atop{R}}\Big)$ quadrants, and so we have
\begin{align*}
\mathbb{P}_{n}(J_1<n,J_2<n,I_1\geq n-R)\leq 2^{R}\Big({{n}\atop{R}}\Big)e^{-(n-1)\tilde{c}}\frac{\int\limits_{\mathcal{Q}_1}f_n(\phi)d\phi}{\int\limits_{\mathcal{Q}_1}f_n(\phi)d\phi}\leq e^{-\Omega(n)}
\end{align*}
completing the proof of the Proposition.

\end{itemize}

\section{Proof of Lemma \ref{l31} and Lemma \ref{lem:G2}}\label{sec:four4}

We first prove an initial estimate, which will be used to verify Lemma  \ref{l31}.
\begin{lem}\label{l30}
There exists finite positive constants $M_1,M_2$ free of $n$ such that
\begin{align}
\label{es1}&\mathbb{P}_n\left(\sum\limits_{i=1}^n(\phi_i-\bar{\phi})^2>M_1\right)\leq e^{-\Omega(n)},\\
\label{es2}&\mathbb{P}_n\left(|\bar{\phi}|>M_2n^{-1/4}\right)\le e^{-\Omega(n)}.
\end{align}
\end{lem}

\begin{proof}
\begin{enumerate}
\item[(a)]
Expanding $q$ around  $\bar{\phi}$ we have
 $$q\left(\frac{\phi_i+\phi_j}{2}\right)=q(\bar{\phi})+\frac{1}{2}(\phi_i+\phi_j-2\bar{\phi})q'(\bar{\phi})+\frac{1}{8}(\phi_i+\phi_j-2\bar{\phi})^2q''(\xi_{ij}).$$ 
 Since $0\leq q''(x)=\tanh^2(2\theta x)\leq 1$, using \eqref{eq:domm} gives 
\begin{align}\label{critic1}\frac{n\theta}{4}\sum\limits_{i=1}^n(\phi_i-\bar{\phi})^2+Nq(\bar{\phi})\leq -\log f_n(\phi)\leq \frac{(2n-1)\theta}{4}\sum\limits_{i=1}^n(\phi_i-\bar{\phi})^2+Nq(\bar{\phi}).
\end{align} 
Let  ${\bf O}_n$ be an orthogonal matrix with first row equal to $\frac{1}{\sqrt{n}}{\bf 1}$. Changing variables to $\psi={\bf O}_n\phi$ gives 
\begin{align}\label{eq:before}
\psi_1=\sqrt{n}\bar{\phi},\quad  \sum\limits_{i=2}^n\psi_i^2=\sum\limits_{i=1}^n(\phi_i-\bar{\phi})^2,
\end{align}
and so \eqref{critic1} becomes
\begin{align*}
\frac{n\theta}{4}\sum_{i=2}^n\psi_i^2+Nq\Big(\frac{1}{\sqrt{n}}\psi_1\Big)\leq -\log f_n(\phi)\leq \frac{n\theta}{2}\sum_{i=2}^n\psi_i^2+Nq\Big(\frac{1}{\sqrt{n}}\psi_1\Big).
\end{align*}
This gives
 \begin{align*}
 \mathbb{P}_{n}\left(\sum\limits_{i=2}^n(\phi_i-\bar{\phi})^2>M_1\right)= &\frac{\int_{\sum_{i=1}^n(\phi_i-\bar{\phi})^2>M_1} f_n(\phi)d\phi}{\int_{\mathbb{R}^n}f_n(\phi)d\phi}\\
 \leq &\frac{\int_{\sum_{i=2}^n\psi_i^2>M_1}e^{-\frac{n\theta}{4}\sum\limits_{i=2}^n\psi_i^2-Nq\Big(\frac{1}{\sqrt{n}}\psi_1\Big)}d\psi}{\int_{\mathbb{R}^n}e^{-\frac{n\theta}{2}\sum\limits_{i=2}^n\psi_i^2-Nq\Big(\frac{1}{\sqrt{n}}\psi_1\Big)}d\psi}
 \leq2^{n/2}\mathbb{P}_n\Big(\chi^2_{n-1}>\frac{nM_1\theta}{2}\Big)
  \end{align*}
Using standard tail bounds for a $\chi^2$ random variable, choosing $M_1$ large enough the above term is $e^{-\Omega(n)}$, completing the proof of part (a).

\item[(b)]
Define a function $r:[-3,3]\mapsto \R$ by 
\begin{align*}
r(t)=&\frac{q(t)}{t^4} \text{ if }t\neq 0,\\
=&\frac{q''''(0)}{4!}=\frac{1}{8}\text{ if }t=0,
\end{align*}
and note that $r$ is continuous and strictly positive on $[-3,3]$. Thus we have
\begin{align}\label{critic2}
0<\lambda_2:=\inf_{|t|\le 3}r(t)\leq \sup_{|t|\le 1}r(t)=:\lambda_1<\infty,
\end{align}
which gives
\begin{align*}
\P_n(|\bar{\phi}|>M_2n^{-1/4})\le \P_n(|\bar{\phi}|>3)+\P_n\Big( |\bar{\phi}|>M_2n^{-1/4}\Big||\bar{\phi}|\le 3\Big)
\end{align*}
where the first term is $e^{-\Omega(n)}$ by \eqref{eq:crude}. Turning to deal with the second term, note that using \eqref{critic1}) along with \eqref{critic2} we get
\begin{align*}
\P_n( |\bar{\phi}|\geq M_2n^{-1/4}\Big||\bar{\phi}|\le 3)\leq &\frac{\int_{M_2n^{-1/4}\leq |\bar{\phi}|\le 3}e^{-\frac{n\theta}{4}\sum_{i=1}^n(\phi_i-\bar{\phi})^2-Nq(\bar{\phi})}d\phi}{\int_{|\bar{\phi}|\le 3} e^{- \frac{(2n-1)\theta}{4}\sum_{i=1}^n(\phi_i-\bar{\phi})^2-Nq(\bar{\phi})}d\phi}\\
=&\frac{\int_{M_2 n^{1/4}\leq |\psi_1|\leq 3\sqrt{n}}e^{-\frac{n\theta}{4}\sum\limits_{i=2}^n\psi_i^2-Nq(\frac{1}{\sqrt{n}}\psi_1)}d\psi}{\int_{|\psi_1|\le 3\sqrt{n}} e^{- \frac{(2n-1)\theta}{4}\sum_{i=2}^n\psi_i^2-Nq(\frac{1}{\sqrt{n}}\psi_1)}d\psi}\\
=&\sqrt{\left(\frac{2n-1}{n}\right)^{n-1}}\frac{\int_{M_2 n^{1/4}\leq |\psi_1|\leq 3\sqrt{n}}e^{-Nq(\frac{1}{\sqrt{n}}\psi_1)}d\psi_1}{\int_{|\psi_1|\leq 3\sqrt{n}}e^{-Nq(\frac{1}{\sqrt{n}}\psi_1)}d\psi_1}\\
\leq &\sqrt{2^{n-1}}\frac{\int_{|\psi_1|\geq  M_2n^{1/4}}e^{-\lambda_2(1-1/n)\psi_1^4/2}d\psi_1}{\int_{|\psi_1|\leq 3\sqrt{n}}e^{-\lambda_1(1-1/n)\psi_1^4/2}d\psi_1}\\
\leq &\sqrt{2^{n-1}}\frac{\int_{|\psi_1|\geq  M_2n^{1/4}}e^{-2^2\lambda_2\psi_1^4/2}d\psi_1}{\int_{|\psi_1|\leq 3\sqrt{n}}e^{-\lambda_1\psi_1^4}d\psi_1/2}.\\
\end{align*}
The denominator of the above ratio bounded away from $0$, whereas the numerator is $e^{-M_2^4\Omega(n)}$, and so if $M_2$ is  large enough then the above term is $e^{-\Omega(n)}$, as desired.
\end{enumerate}
\end{proof}
\subsection{ {Proof of Lemma \ref{l31}} part (a)}

Setting $\bar{\phi}_1:=(n-1)^{-1}\sum\limits_{i=2}^n\phi_i$, note that \begin{align}\label{eq:crude2}
|\bar{\phi}-\bar{\phi}_1|\le \frac{|\sum_{i=2}^n\phi_i|}{n(n-1)}+\frac{|\phi_1|}{n}\le \frac{6}{n},
\end{align}
where the last inequality holds on the set $\max_{1\le i\le n}|\phi_i|\le 3$. This gives
\begin{align*}
\E |\bar{\phi}-\bar{\phi}_1|^\ell \lesssim_\ell \frac{1}{n^\ell}+\E\Big[ |\bar{\phi}-\bar{\phi}_1|^\ell, 1\{\max_{1\le i\le n}|\phi_i|>3)\Big]\lesssim_\ell \frac{1}{n^\ell},
\end{align*}
where the last inequality uses \eqref{eq:crude}. It thus suffices to show that
\begin{align}\label{es2.5}
\E|\phi_1-\bar{\phi}_1|^l\lesssim_\ell n^{-l/2}.
\end{align}
To this effect, use \eqref{eq:crude2} along with part (a) of Lemma \ref{l30} to get
the existence of $M_3<\infty$ such that
\begin{align}\label{es3}
\P_n(|\bar{\phi}_1|>M_3n^{-1/4})\le e^{-\Omega(n)}.
\end{align}
Now, for $2\leq i\leq n$ expanding $q\Big(\frac{\phi_1+\phi_i}{2}\Big)$ around $\bar{\phi}_1$ we get
\[q\Big(\frac{\phi_1+\phi_i}{2}\Big)=q(\bar{\phi}_1)+\frac{1}{2}(\phi_1+\phi_i-2\bar{\phi}_1)q'(\bar{\phi}_1)+\frac{1}{8}(\phi_1+\phi_i-2\bar{\phi}_1)^2q''(\xi_{ij})\] with $0\leq q''(x)\le 1$ as before. This along with \eqref{eq:domm} gives that 
\begin{align*}
&q(\bar{\phi}_1)+\frac{1}{2}(\phi_1+\phi_i-2\bar{\phi}_1)q'(\bar{\phi}_1)+\frac{1}{8}(\phi_1-\phi_i)^2
\le p(\phi_1,\phi_i)\\
\le &q(\bar{\phi}_1)+\frac{1}{2}(\phi_1+\phi_i-2\bar{\phi}_1)q'(\bar{\phi}_1)+ \frac{1}{8}(\phi_1-\phi_i)^2+\frac{1}{8}(\phi_1+\phi_i-2\bar{\phi}_1)^2.
\end{align*}
On adding this over $2\le i\le n$ gives
\begin{align}
\notag&(n-1)q(\bar{\phi}_1)+(n-1)(\phi_1-\bar{\phi}_1)q'(\bar{\phi}_1)+\frac{n-1}{8}(\phi_1-\bar{\phi}_1)^2+\frac{1}{8}\sum\limits_{i=2}^n(\bar{\phi}_1-\phi_i)^2
\notag\leq \sum\limits_{i=2}^np(\phi_1,\phi_i)\\
\label{es4}\leq&(n-1)q(\bar{\phi}_1)+(n-1)(\phi_1-\bar{\phi}_1)q'(\bar{\phi}_1)+ \frac{n-1}{4}(\phi_1-\bar{\phi}_1)^2+\frac{1}{4}\sum\limits_{i=2}^n(\bar{\phi}_1-\phi_i)^2.
\end{align}
Setting $$A_n(M):=\Big\{\phi\in \R^n: \sum\limits_{i=2}^n(\phi_i-\bar{\phi}_1)^2\leq M,\quad |\bar{\phi}_1|\leq Mn^{-1/4} \Big\},$$
we can write
\begin{align}
\begin{split}
\label{eq:sol_1}
\E|\phi_1-\bar{\phi}_1|^l=\E\left[|\phi_1-\bar{\phi}_1|^l,1\{A_n(M)\}\right]
+\E\left[|\phi_1-\bar{\phi}_1|^l,1\{A_n(M)^c\}\right].
\end{split}
\end{align}
Using Cauchy-Schwarz inequality, the second term in RHS of \eqref{eq:sol_1} can be bounded as follows:
\begin{align} 
\begin{split}\label{eq:sol_2}
\E\left[|\phi_1-\bar{\phi}_1|^l,1\{A_n(M)^c\}\right]
\le\sqrt{\E|\phi_1-\bar{\phi}_1|^{2l}}\sqrt{\P_n(A_{n}(M))}\le e^{-\Omega(n)},
\end{split}
\end{align} 
where the last equality uses Lemma \ref{l30} to conclude that $\P_n(A_n(M)^c)\le e^{-\Omega(n)}$ for some $M$ fixed. For the first term in the RHS of \eqref{eq:sol_1} we have
\begin{align*}
\E\left[|\phi_1-\bar{\phi}_1|^l,1\{A_n(M)\}\right]
=\E\left[\E\Big(|\phi_1-\bar{\phi}_1|^l\Big|\phi_i,2\leq i\leq n\Big)1\{A_n(M)\}\right],
\end{align*}
with
\begin{align}
\notag&\E\Big(|\phi_1-\bar{\phi}_1|^l\Big|\phi_i,2\leq i\leq n\Big)
=\frac{\int\limits_{\R}|\phi_1-\bar{\phi}_1|^l\prod\limits_{i=2}^ne^{-p(\phi_1,\phi_i)}d\phi_1}{\int\limits_{\R}\prod\limits_{i=2}^ne^{-p(\phi_1,\phi_i)}d\phi_1}\\
\label{last}\le &e^{\frac{1}{8}\sum\limits_{i=2}^n(\phi_i-\bar{\phi}_1)^2}\frac{\int\limits_\R|\phi_1-\bar{\phi}_1|^le^{-(n-1)(\phi_1-\bar{\phi}_1)q'(\bar{\phi}_1)-\frac{n-1}{8}(\phi_1-\bar{\phi}_1)^2}d\phi_1}{\int\limits_\R e^{-(n-1)(\phi_1-\bar{\phi}_1)q'(\bar{\phi}_1)-\frac{n-1}{4}(\phi_1-\bar{\phi}_1)^2}d\phi_1},
\end{align}
where we have used \eqref{es4} in the last step. 

 Focusing on the ratio in \eqref{last}, the numerator can be simplified as
\begin{align*}
\int\limits_\R|\phi_1-\bar{\phi}_1|^le^{-(n-1)(\phi_1-\bar{\phi}_1)q'(\bar{\phi}_1)-\frac{n-1}{8}(\phi_1-\bar{\phi}_1)^2}d\phi_1=&\int\limits_{\R}|z|^le^{-\frac{(n-1)}{2.4}[2.4zq'(\bar{\phi}_1)+z^2]}dz\\
=&e^{2(n-1)q'(\bar{\phi}_1)^2}\sqrt{\frac{8\pi}{n-1}}\mathbb{E}|\widetilde{W}|^l
\end{align*}
where $\widetilde{W}\sim N\Big(-4q'(\bar{\phi}_1),\frac{4}{n-1}\Big)$.  Thus the numerator in \eqref{last} is bounded (upto constants depending on $\ell$) by
\[n^{-1/2}e^{2(n-1)q'(2\bar{\phi})^2}\bigg[|q'(\bar{\phi}_1)|^l+\frac{1}{n^{l/2}}\bigg].\]
By a similar calculation, the denominator is lower bounded by (up to universal constants)
$n^{-1/2}e^{(n-1)q'(\bar{\phi}_1)^2}.$ Using \eqref{last} then gives
\begin{align}\label{eq:last}
\E\Big(|\phi_1-\bar{\phi}_1|^l\Big|\phi_i,2\leq i\leq n\Big)\lesssim_\ell e^{nq'(\bar{\phi}_1)^2}\bigg[|q'(\bar{\phi}_1)|^l+n^{-l/2}\bigg].
\end{align}
Finally, using arguments as before, 
  there exists finite positive constants $\lambda_1\ge \lambda_2$ such that for all $x\in [-1,1]$ we have $ \lambda_2 |x|^3\leq q'(x)\leq \lambda_1 |x|^3.$ Since $|\bar{\phi}_1|\leq M_3n^{-1/4}$ on the conditioned set, for all $n$ large we have
$$e^{nq'(\bar{\phi}_1)^2}[|q'(\bar{\phi}_1)|^l+n^{-l/2}|\lesssim_\ell [n^{-3l/4}+n^{-l/2}]\lesssim_\ell n^{-l/2}.$$
Combining \eqref{eq:last} with the above display gives that on the
\begin{align*}
\E\Big(|\phi_1-\bar{\phi}_1|^l\Big|\phi_i,2\leq i\leq n\Big)\lesssim n^{-\ell/2}.
\end{align*}
It follows from this
 $\E|{\phi}_1-\bar{\phi}_1|\lesssim n^{-l/2}$, proving \eqref{es2.5}, and hence completing the proof of part (a).

\subsection{ {Proof of Lemma \ref{l31} part (b)}}

Using \eqref{eq:phi_z}, it suffices to show that $$\limsup\limits_{n\rightarrow\infty}\mathbb{E}|S(Y)|^4<\infty,\quad S(Y):=n^{-3/2}\sum_{e\in \mathcal{E}} Y_e.$$
Produce the usual exchangeable pair $(Y,Y')$ as in Lemma \ref{pair}, and note that
\begin{align}\label{curie1}
\mathbb{E}(S(Y)-S(Y')|Y=y)=\frac{1}{Nn\sqrt{n}}\sum\limits_{e\in E}\Big[y_e-\tanh( t_e)\Big]\end{align}
For $e\in E$, setting $\bar{t}:=\frac{1}{N}\sum\limits_{e\in E}t_e$ and expanding $\tanh$ by a Taylor's series gives
\begin{align}\label{curie2}
\tanh( t_e)=\tanh( \bar{t})+(t_e-\bar{t})\text{ sech}^2( \bar{t})
-(t_e-\bar{t})^2\text{sech}^2(\bar{t})\tanh( \bar{t})+\frac{f(\xi_e)(t_e-\bar{t})^3}{3!}
\end{align}
for some $\xi_e$, with $|f(x)|:=|\tanh'''(x)|\leq 1$. Summing over $e\in E$ we have \begin{align*}
\sum\limits_{e\in E}\tanh(t_e)=&N\tanh(\bar{t})-\text{ sech}^2(\bar{t})\tanh(\bar{t})\sum_{e\in E}(t_e-\bar{t})^2
+\frac{1}{3!}\sum_{e\in E}(t_e-\bar{t})^3f(\xi_e)
\end{align*}
Using the fact that $|\tanh(x)-x+x^3/3|\le {16}|x|^5/{5!}$ for all $x\in \R$, we have 
\begin{align}\label{curie3}
\Big|\tanh(\bar{t})-\bar{t}+\frac{(\bar{t})^3}{3}\Big|\leq \frac{16}{5!}|\bar{t}|^5
\end{align}
Plugging in the estimates \eqref{curie2} and \eqref{curie3} in \eqref{curie1} gives
\begin{align}\label{curie4}
\notag&\mathbb{E}(S(Y)-S(Y')|Y=y)\\
=&\frac{(\bar{t})^3}{3n\sqrt{n}}+\frac{1}{Nn\sqrt{n}}\text{ sech}^2(\bar{t})\tanh(\bar{t})\sum\limits_{e\in E}(t_e-\bar{t})^2+A_1+A_2+A_3,
\end{align}
with 
\begin{align}\label{eq:first_bound}
\begin{split}
|A_1|=&\Big|\frac{1}{Nn\sqrt{n}}\sum_{e\in E}(y_e-\bar{t})\Big|=\frac{|\sum_{e\in E}y_e|}{N(n-1) n\sqrt{n}}\lesssim \frac{|S(Y)|}{n^3},\\
|A_2|\le & \frac{2}{Nn\sqrt{n}}\sum_{e\in E}|t_e-\bar{t}|^3\lesssim \frac{1}{n^{3.5}}\sum_{e\in E}|t_e-\bar{t}|^3,\\
|A_3|\le &\frac{16|\bar{t}|^5}{5! n\sqrt{n}}\le \Big(\frac{2^9}{5!}+o(1)\Big)\frac{|S(Y)|^5}{n^4}.
\end{split}
\end{align}
Using \eqref{curie4} gives
{\small \begin{align}\label{curie5}
\begin{split}
\mathbb{E}\Big[(S(Y)-S(Y'))S(Y)\Big]=&\mathbb{E}\Big[\mathbb{E}(S(Y)-S(Y')|Y)S(Y)\Big]\\
=&\frac{1}{3n\sqrt{n}}\mathbb{E}\Big[S(Y)\bar{t}^3\Big]+\mathbb{E}\Big[(A_1+A_2+A_3)S(Y)\Big]\\
+&\frac{1}{Nn\sqrt{n}}\mathbb{E}\Big[\text{ sech}^2(\bar{t})\tanh(\bar{t})\sum\limits_{e\in E}(t_e-\bar{t})^2S(Y)\Big].
\end{split}
\end{align}}
 We now bound each term on the RHS of \eqref{curie5}. To begin, use  \eqref{eq:first_bound} to note that
 \begin{align}\label{est1}
 \begin{split}
&\left| \mathbb{E}\Big[A_1S(Y)\Big]\right|\lesssim n^{-3}\mathbb{E}S(Y)^2,\\
&\left | \mathbb{E}\Big[A_3S(Y)\Big]\right|\le \Big(\frac{2^9}{5!}+o(1)\Big)n^{-4}\mathbb{E}S(Y)^6\le \Big(\frac{2^7}{5!}+o(1)\Big)n^{-3}\mathbb{E} S(Y)^4,
\end{split}
\end{align}
where the last inequality uses the trivial bound $|S(Y)|\le \frac{\sqrt{n}}{2}$.
 Next we have
\begin{align}\label{est2}
\begin{split}
\left| \mathbb{E}\Big[A_2S(Y)\Big]\right|\lesssim& n^{-3.5}\sum\limits_{e\in E}\mathbb{E}\Big[|S(Y)||t_e-\bar{t}|^3\Big]\\
\le &n^{-3.5}\sqrt{\mathbb{E}S(Y)^2}\sum_{e\in E}\sqrt{\mathbb{E}(t_e-\bar{t})^6}
\lesssim n^{-3}\sqrt{\mathbb{E}S(Y)^2}.
\end{split}
 \end{align}
Here the last inequality follows from the fact that 
 \begin{align}\label{eq:l_moment}
 \mathbb{E}|t_e-\bar{t}|^\ell \lesssim n^{-\ell }\mathbb{E}|k_1-\bar{k}|^\ell\lesssim n^{-\ell} \E |\phi_1-\bar{\phi}|^\ell \lesssim n^{-\frac{\ell}{2}},
 \end{align} where the last estimate uses part (b) of Lemma \ref{l31}. 
 Proceeding to estimate the final term in the RHS of \eqref{curie5} we have
 \begin{align}
 \begin{split}\label{est3}
 \frac{1}{Nn\sqrt{n}}\mathbb{E}\Big[\text{ sech}^2(\bar{t})\tanh(\bar{t})\sum\limits_{e\in E}(t_e-\bar{t})^2S(Y)\Big]
\lesssim&n^{-4}\sum\limits_{e\in E}\mathbb{E}S(Y)^2|t_e-\bar{t}|^2\\
\lesssim&n^{-4}\sqrt{\mathbb{E}S(Y)^4}\sum_{e\in E}\sqrt{\mathbb{E}|t_e-\bar{t}|^4}\\
 \lesssim &n^{-3}\sqrt{\mathbb{E}S(Y)^4},
 \end{split}
 \end{align}
 where the last estimate again uses \eqref{eq:l_moment}.
  Since \[\frac{1}{3n\sqrt{n}}\mathbb{E}[S(Y)\bar{t}^3]=\Big(\frac{8}{3}+o(1)\Big)\frac{\mathbb{E}S(Y)^4}{n^3},\]  using the estimates obtained in  \eqref{est1}, \eqref{est2} and \eqref{est3} along with \eqref{curie5} we have
  \begin{align*}
 &\left| n^3\mathbb{E}(S(Y)-S(Y'))S(Y)-\Big(\frac{8}{3}+o(1)\Big)\mathbb{E}S(Y)^4\right|-\Big(\frac{2^7}{5!}+o(1)\Big)\E S(Y)^4\\
 \lesssim&
 \sqrt{\mathbb{E}S(Y)^4}+\sqrt{\mathbb{E}S(Y)^2}\lesssim \sqrt{\E S(Y)^4},
  \end{align*}
  where the last bound uses H\"older's inequality. Noting that $\frac{8}{3}>\frac{2^7}{5!}$, this gives
\begin{align}\label{est4}
\E S(Y)^4\lesssim n^3 \left|\E \Big[\Big(S(Y)-S(Y'))S(Y)\Big]\right|+\sqrt{\E S(Y)^4}.
  \end{align}
Finally, by exchangeability of $(Y,Y')$ we have
\begin{align*}
\mathbb{E}\Big[(S(Y)-S(Y'))S(Y)\Big]=\mathbb{E}\Big[(S(Y')-S(Y))S(Y')\Big]=\frac{1}{2}\mathbb{E}\Big[(S(Y)-S(Y'))^2\Big]
\lesssim n^{-3},
\end{align*}
which along with \eqref{est4} gives the existence of a constant $C$ free of $n$ such that
\begin{align*}
\mathbb{E}S(Y)^4\le C\Big[1+\sqrt{\E S(Y)^4}\Big].
\end{align*}
The desired conclusion follows from this.

\subsection{Proof of Lemma \ref{lem:G2}}
Let $\zeta$ be a random variable with density proportional to $e^{-\frac{\zeta^2}{2}-\frac{\zeta^4}{24}}$, as in the statement of the lemma. Let $W_1,\ldots,W_n\stackrel{i.i.d.}{\sim}N(0,1)$ independent of $\zeta$. Then we claim that 
\begin{align}
\label{eq:rep_c}
(\phi_1,\ldots,\phi_n)\stackrel{D}{=}\zeta+\frac{2}{\sqrt{n-1}}(W_1-\bar{W},\ldots,W_n-\bar{W}).
\end{align} 
Indeed, to verify this, let  ${\bf O}_n$ is an orthogonal matrix with first row equal to $n^{-1/2}{\bf 1}$, as in the proof of Lemma \ref{l30}. Thus, as in \eqref{eq:before}, we have 
\begin{align}\label{eq:before2}
\psi_1=\sqrt{n}\bar{\phi},\quad \sum_{i=2}^n\psi_i^2=\sum_{i=1}^n(\phi_i-\bar{\phi})^2.
\end{align}
Also, the joint distribution of $\psi$ is proportional to $\overline{g3}(.)$, where
\begin{align}\label{g_bar}
-\log \overline{g}_{3n}(\psi)=\frac{(n-1)}{8}\sum\limits_{i=2}^n\psi_i^2+\frac{1}{24}\psi_1^4+\frac{1}{2}\psi_1^2.    
\end{align}
In particular this means that $\{\psi_i\}_{1\le i\le n}$ are mutually independent, with $\psi_1\stackrel{D}{=}\zeta$, and $\psi_i\sim N(0,\frac{4}{n-1})$ for $2\le i\le n$. Since $(\phi_1-\bar{\phi},\ldots,\phi_n-\bar{\phi})$ is a linear transformation of $(\psi_2,\ldots,\psi_n)$, it follows that $(\phi_1-\bar{\phi},\ldots,\phi_n-\bar{\phi})$ has a multivariate (singular) Gaussian distribution with mean vector ${\bf 0}$ and covariance matrix $\Sigma$, where
\begin{align*}
&\Sigma_{ii}={\sf Var}(\phi_i-\bar{\phi})=\frac{n-1}{n}{\sf Var}(\psi_2)=\frac{4}{n},\\
&\Sigma_{ij}=-\frac{1}{n-1}\Sigma_{ii}=-\frac{4}{n(n-1)}.
\end{align*}
By matching the covariance, it follows that $$(\phi_1-\bar{\phi},\ldots,\phi_n-\bar{\phi})\stackrel{D}{=}\frac{2}{\sqrt{n-1}}(W_1-\bar{W},\ldots,W_n-\bar{W}),$$
where $\{W_i\}_{1\le i\le n}\stackrel{i.i.d.}{\sim}N(0,1)$, as in the statement of the lemma. We have thus verified \eqref{eq:rep_c}, which we now use to verify all parts of the lemma. Before proceeding, note that part (c) is immediate from the above calculations.

\begin{enumerate}
\item[(a)]
\eqref{eq:fg2}  follows on using \eqref{eq:before2} to note that 
$n\E\bar{\phi}^2=\E \zeta^2$, and $n\E(\phi_1-\bar{\phi})^2\lesssim n\E \psi_2^2\lesssim 1$.

\item[(b)]
\eqref{eq:g21} follows on using \eqref{eq:before2} to note that
$\sum_{i=1}^n(\phi_i-\bar{\phi})^2=\sum_{i=2}^n\psi_i^2$, which converges to $4$ in probability.
For verifying \eqref{eq:g22} and \eqref{eq:g23}, for every positive integer $r$ using \eqref{eq:rep_c}  we have
\begin{align}\label{eq:first_note}
\sum_{i<j}(\phi_i+\phi_j-2\bar{\phi})^r\stackrel{D}{=}\Big(\frac{2}{\sqrt{n-1}}\Big)^r\sum_{i<j}(W_i+W_j-2\bar{W})^r.
\end{align}
The desired conclusion then follows on noting the following two limits:
\begin{align*}
\frac{1}{n^2}\sum_{i<j}(W_i+W_j-2\bar{W})^3\stackrel{P}{\to}0,\\
\frac{1}{n^2}\sum_{i<j}(W_i+W_j-2\bar{W})^4\stackrel{P}{\to}6.
\end{align*}

\item[(d)]

This is immediate on using \eqref{eq:before2} to note that $$\sum_{i=1}^n(\phi_i-\bar{\phi})^2=\sum_{i=2}^n\psi_i^2\stackrel{D}{=}\frac{4}{n-1}\chi_{n-1}^2,$$
and then using a CLT for $\chi_{n-1}^2$.

\item[(e)]
This follows on using \eqref{eq:rep_c} to note that
$$\sum_{i=1}^nc_n(i)\phi_i\stackrel{D}{=}\frac{2}{\sqrt{n-1}}\sum_{i=1}^nc_n(i)(W_i-\bar{W})=\frac{2}{\sqrt{n-1}}\sum_{i=1}^nc_n(i)W_i\sim N\Big(0,\frac{4}{n-1}\sum_{i=1}^nc_n(i)^2\Big).$$

\item[(f)]
This follows on using \eqref{eq:rep_c} to note that
$$\max_{i,j\in [n]:|i-j|\le n\delta}|S_i(\phi)-S_j(\phi)|\stackrel{D}{=}\frac{2}{\sqrt{n-1}}\max_{i,j\in [n]:|i-j|\le n\delta}|\widetilde{S}_i({\bf W})-\widetilde{S}_j({\bf W})|,$$
where ${\bf W}:=(W_1,\ldots,W_n)$, and $\widetilde{S}_i(\phi)=\sum_{j=1}^iW_j$ as in the proof of Lemma \ref{lem:G1}. The desired conclusion follows as before on using sample path tightness for partial sums of i.i.d.~random variables.

\end{enumerate}

\end{appendix}

\end{document}